\documentclass[12pt, reqno]{amsart}
\usepackage[utf8]{inputenc}
\usepackage[english]{babel}
\usepackage{comment}
\usepackage{blkarray}
\usepackage{fullpage}              
\usepackage{tikz}    
\usetikzlibrary{arrows.meta,arrows,patterns}
\usetikzlibrary{shapes,decorations,decorations.pathreplacing,arrows,calc,arrows.meta,fit,positioning}
\usetikzlibrary{arrows.meta,automata,positioning,quotes}
\usepackage[pdftex]{hyperref}
\hypersetup{
    colorlinks=true,
    linkcolor={blue},
    citecolor={magenta},
    urlcolor={blue}
}
\usepackage[nameinlink,capitalise]{cleveref}
\usepackage{mathtools}

\DeclareMathOperator{\de}{de}

\numberwithin{equation}{section}
\theoremstyle{plain}
\newtheorem{theorem}{Theorem}[section]

\newtheorem{conjecture}{Conjecture}[section]
\newtheorem{corollary}[theorem]{Corollary}
\newtheorem{lemma}[theorem]{Lemma}
\newtheorem{proposition}[theorem]{Proposition}
\newtheorem{problem}[theorem]{Problem}

\theoremstyle{definition}
\newtheorem{remark}[theorem]{Remark}
\newtheorem{examplex}[theorem]{Example}
\newenvironment{example}
  {\pushQED{\qed}\examplex}
  {\popQED\endexamplex}

\newtheorem{definition}[theorem]{Definition}

\newcommand{\cc}{\mathbb{C}}
\newcommand{\rr}{\mathbb{R}}

\newcommand{\cm}{\mathcal{M}}

\newcommand{\cg}{\mathcal{G}}

\newcommand{\cv}{\mathcal{V}}

\DeclareMathOperator{\image}{image}

\DeclareMathOperator{\pa}{pa}
\DeclareMathOperator{\an}{an}
\DeclareMathOperator{\nd}{nd}

\DeclareMathOperator{\glob}{global}
\newcommand{\indep}{\perp\!\!\!\perp}
\newcommand{\pd}{\mathrm{PD}}

\newcommand{\gmci}{\cm_{\cg, i \indep j|K}}
\newcommand{\vgmci}{V_{\cg, i \indep j|K}}

\definecolor{benpurple}{RGB}{180, 0, 240}

\author{Mathias Drton}
\address{Technische Universit\"at M\"unchen, 85748 Garching b. München, Boltzmannstr. 3.,  Germany}
\email{mathias.drton@tum.de}

\author{Leonard Henckel}
\address{University College Dublin, Belfield, Dublin 4, Ireland}
\email{leonard.henckel@ucd.ie}

\author{Benjamin Hollering}
\address{Technische Universit\"at M\"unchen, 85748 Garching b. München, Boltzmannstr. 3.,  Germany}
\email{benjamin.hollering@tum.de}

\author{Pratik Misra}
\address{Technische Universit\"at M\"unchen, 85748 Garching b. München, Boltzmannstr. 3.,  Germany}
\email{pratik.misra@tum.de}

\title{Faithlessness in Gaussian graphical models}

\date{}
\keywords{}
\subjclass{}

\begin{document}

\begin{abstract}
The implication problem for conditional independence  (CI) asks whether the fact that a probability distribution obeys a given finite set of CI relations implies that a further CI statement also holds in this distribution.  This problem has a long and fascinating history, cumulating in positive results about implications now known as the semigraphoid axioms as well as impossibility results about a general finite characterization of CI implications.  Motivated by violation of faithfulness assumptions in causal discovery, we study the implication problem in the special setting where the CI relations are obtained from a directed acyclic graphical (DAG) model along with one additional CI statement.
%
Focusing on the Gaussian case, we give a complete characterization of when such an implication is graphical by using algebraic techniques.  
Moreover, prompted by the relevance of strong faithfulness in statistical guarantees for causal discovery algorithms, we give a graphical solution for an approximate CI implication problem, in which we ask whether small values of one additional partial correlation entail small values for yet a further partial correlation.
\end{abstract}

\maketitle

\section{Introduction}

Conditional independence is an important tool in probabilistic and causal reasoning. The set of conditional independencies that hold in a large random vector $X$ can, for example, be used to efficiently store and compute complex conditional probabilities in $X$; a task that would otherwise be infeasible in large dimensions \cite{koller2009probabilistic}. In order to effectively use conditional independence as a tool, it is important to understand, given a set of conditional independence statements which other independence statements have to hold in any distribution that satisfies the former. This problem is known as the implication problem and has been extensively studied. In particular, it has been shown that for an important class of independence statements, called a recursive basis, the implication problem can be solved by applying four rules called the semi-graphoid axioms \cite{geiger1993logical,geiger1990identifying,studeny1997semigraphoids}. Recursive bases are used to define probabilistic graphical models such as directed acyclic graphical models (DAGs), where all conditional independence statements implied by the semi-graphoid axioms can be found by applying the d-separation criterion \cite{geiger1990logic,geiger1990identifying}. For general collections of conditional independence statements, however, it has been shown that no finite axiomatization to solve the implication problem exists \cite{studeny1992conditional}. As a result, most research on probabilistic and causal reasoning has focused on probabilistic graphical models, where the implication problem can be solved in polynomial time \cite{koller2009probabilistic,pearl2009causality,spirtes2000causation}. 

In order to understand how strong this restriction to probabilistic graphical models is, it is important to study how stable they are as a class of independence models and whether there are intermediate models between probabilistic graphical models and the general case in terms of complexity. In this paper, we aim to study this question from an algebraic perspective by investigating how a DAG model for the conditional independencies of a Gaussian random vector behaves if we add an additional conditional independence statement which is not implied by d-separation. In particular we investigate under what conditions on the added conditional independence statement, the model remains a DAG-model and when the implication problem remains solvable with the semi-graphoid axioms. 

The class of DAG-models in which at least one additional conditional independence statements holds is also interesting from a causal perspective. One standard approach to causal discovery is to learn the causal graph from conditional independence statements \cite{sadeghi2022conditions, spirtes2001anytime,spirtes2000causation}. Most algorithms that use this approach are only consistent for causal models where no conditional independence statements other than those implied by d-separation in the true causal DAG hold; an assumption called faithfulness. Studying the implication problem in a DAG model with one additional conditional independence statement, therefore corresponds to studying the algebraic geometry of a faithfulness violation; a problem we call faithlessness propagation.

Both the faithfulness assumption and the implication problem have been criticized for only considering exact conditional independence. In finite samples, a very small dependence cannot be distinguished statistically from exact independence. In response, a literature has evolved studying a stronger version of faithfulness known as $\lambda$-faithfulness \cite{uhler2013geometry,zhang2002strong}. Similarly, there exists a literature studying the implications of approximate conditional independence statements \cite{kenig21approximate,kenig2022integrity}. The latter has established that while an approximate implication implies exact implication the reverse is not necessarily the case. As such, the question of how a $\lambda$-faithfulness violation propagates in a DAG model is distinct from the problem of faithlessness propagation. The problem of $\lambda$-faithlessness is also of independent interest for the following reason: conditional independence statements in causal models can be used to identify more efficient causal effect estimators \cite{henckel2022graphical,henckel2019graphical,rotnitzky2020efficient,runge2021necessary,witte2020efficient}. These results exploit that a conditional independence statement of the form $X \indep A | B$ implies that the mutual information between $X$ and $B$ is larger than that between $X$ and $A$. Therefore, if such a d-separation statement holds in a graph $\cg$ the corresponding mutual information inequality holds for all models compatible with $\cg$. However, d-separation is not a necessary condition for this; $\lambda$-faithlessness propagation on the other hand is. As such a better understanding of $\lambda$-faithlessness propagation may lay the groundwork for stronger efficient causal effect estimation results via a better understanding of information theoretic inequalities in graphical models. 

In this paper, we study the problems of faithlessness and $\lambda$-faithlessness propagation for linear Gaussian DAG-models. We begin with some preliminaries on linear structural equation models in \cref{sec:SEMBackground} that will be needed for the rest of the paper. We state the CI implication problem in \cref{sec: unfaithful} and use algebraic tools to convert the problem into a principal ideal membership problem. Using the results in \cite{positivity2013}, we provide a combinatorial criterion in \cref{sec:graphicalImplication} to check whether we can obtain a graphical implication after the addition of an extra CI statement. In \cref{sec:gaussoid}, we show that the implications obtained from our special set of CI statements can be deduced from the gaussoid axioms for $n=4$, and conjecture that to be true for all $n$. In the end, we give a graphical solution for an approximate CI implication problem in \cref{sec:approximate CI implication}, where the additional partial correlation is arbitrarily small instead of zero.

\section{Preliminaries}
\label{sec:SEMBackground}

In this section we provide background on directed Gaussian graphical models, also known as linear Gaussian structural equation models (SEMs),  which will be one of the primary objects we study throughout this paper. We then describe more general conditional independence models which strictly include graphical models. Throughout this section, we also discuss the algebraic perspective on all of these statistical models. For additional background on graphical models we refer the reader to \cite{LauritzenBook, GMhandbook} and for an in depth look at their algebraic structure we refer the reader to \cite{algstat2018}. Lastly, for a detailed discussion on conditional independence implication problems we highly recommend \cite{TobiasThesis}. We introduce the standard terminology in \cref{appendix: prelims} for completeness.

 Let $\cg = (V, E)$ be a directed acyclic graph (DAG) with node set $V$ and edge set $E$. Let $\epsilon = (\epsilon_i ~|~ i \in V)$ be a vector of independent Gaussian errors with mean $0$ and diagonal covariance matrix $\Omega$, which at times we will identify its vector of positive diagonal entries.  Then a random vector $X$ is distributed according to the linear structural equation model on $\cg$ if it satisfies the recursive structural equation system
\begin{equation*}
    X_j = \sum_{i \in \pa(j)} \lambda_{ij} X_i + \epsilon_i,
\end{equation*}
where $\lambda_{ij}$ are edge weights which are called the \emph{direct causal effect} of $X_i$ on $X_j$ and $\pa(j) = \{i \in V ~|~ i \to j \in E\}$ is the set of \emph{parents} of $j$. This system of recursive equations can be solved explicitly which yields $X = (I - \Lambda)^{-T} \epsilon$ where $\Lambda = (\lambda_{ij})$ is the matrix of edge weights such that $\lambda_{ji} \neq 0$ if and only if $j \to i \in E$. Thus any random vector $X$ which satisfies this system of equations will also be Gaussian with mean 0 and covariance matrix $\Sigma$ given by
\[
\Sigma = (I - \Lambda)^{-T} \Omega (I - \Lambda)^{-1},
\]
where $\Omega$ is the covariance matrix of the vector $\epsilon$ and $\Lambda = (\lambda_{ij})$ is the matrix of edge weights.

\begin{definition}\label{definition:Gaussian DAG models}
Let $\cg=(V,E)$ be a DAG with $n$ nodes.  Identifying centered Gaussian joint distributions with their covariance matrix, we define the \emph{Gaussian linear structural equation model} on $\cg$ as the image of the map
\begin{align*}
\phi_\cg: \rr^E \times (0,\infty)^n &\to \pd_n \\
    (\Lambda,\Omega) &\to (I-\Lambda)^{-T} \Omega (I-\Lambda)^{-1},
\end{align*}
where $\pd_n$ is the cone $n \times n$ positive definite matrices. We denote this model by $\cm_\cg = \image(\phi_\cg)$.  Subsequently, we will refer to $\cm_\cg$ simply as the Gaussian DAG model on $\cg$.
\end{definition}

These models are commonly used throughout statistics and frequently go by other names including Gaussian Bayesian networks or directed Gaussian graphical models. 
The following example illustrates the parameterization of the model described above. 

\begin{figure}
\centering
\begin{tikzpicture}[>=stealth',shorten >=1pt,auto,node distance=0.8cm,scale=.9, transform shape,align=center,minimum size=3em]
\node[state] (w1) at (-2,0) {$1$};
\node[state] (w2) at (0,0) {$2$};
\node[state] (w3) at (2,0) {$4$};
\node[state] (w4) at (4,0) {$3$}; 
\node[state] (w5) at (3,2) {$5$};


\path[->]   (w1) edge    (w2);
\path[->]   (w2) edge    (w3);
\path[<-]   (w3) edge    (w4);
\path[->]   (w3) edge    (w5);
\path[->]   (w4) edge    (w5);

\end{tikzpicture}
\caption{A directed acyclic graph $\cg$.}
\label{figure:running-example}
\end{figure}
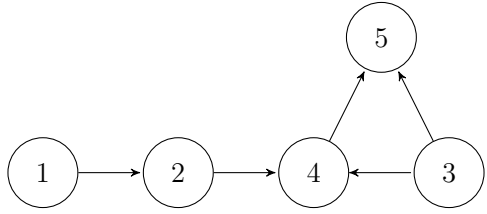

\begin{example}\label{example:structural equations}
Consider the DAG $\cg$ in \cref{figure:running-example}. The recursive structural equations defined by $\cg$ are as follows:
\begin{eqnarray*}
    X_1&=& \epsilon_1, \\
    X_2&=& \lambda_{12}X_1 + \epsilon_2, \\
    X_3&=& \epsilon_3, \\
    X_4&=& \lambda_{24}X_2 + \lambda_{34}X_3 + \epsilon_4 ,\\
    X_5&=& \lambda_{35}X_3 + \lambda_{45}X_4 + \epsilon_5. 
\end{eqnarray*}
The covariance matrix $\Sigma$ for this system can be written as
\[
\Sigma=\begin{bmatrix}
    1 & -\lambda_{12} & 0 & 0 & 0 \\
    0 & 1 & 0 & -\lambda_{24} & 0 \\
    0 & 0 & 1 & -\lambda_{34} & -\lambda_{35} \\
    0 & 0 & 0 & 1 & -\lambda_{45}\\
    0 & 0 & 0 & 0 & 1
\end{bmatrix}^{-T}
\begin{bmatrix}
    \omega_1 & 0 & 0 & 0 & 0 \\
    0 & \omega_2 & 0 & 0 & 0 \\
    0 & 0 & \omega_3 & 0 & 0 \\
    0 & 0 & 0 & \omega_4 & 0 \\
    0 & 0 & 0 & 0 & \omega_5
\end{bmatrix}
\begin{bmatrix}
    1 & -\lambda_{12} & 0 & 0 & 0 \\
    0 & 1 & 0 & -\lambda_{24} & 0 \\
    0 & 0 & 1 & -\lambda_{34} & -\lambda_{35} \\
    0 & 0 & 0 & 1 & -\lambda_{45}\\
    0 & 0 & 0 & 0 & 1
    \end{bmatrix}^{-1}.
\]
Expanding this product allows us to write the covariances in the following way:
\begin{eqnarray*}
    \sigma_{11}&=& \omega_1, \\
    \sigma_{12}&=& \omega_1 \lambda_{12}, \\
    \sigma_{13}&=& 0, \\
    \sigma_{14}&=& \omega_1\lambda_{12}\lambda_{24}, \\
    &&\vdots \\
    \sigma_{55}&=& \omega_{1} \lambda_{12}^2\lambda_{24}^2\lambda_{45}^2 + \omega_2\lambda_{24}^2\lambda_{45}^2 + \omega_3\lambda_{34}^2\lambda_{45}^2 \\
    &&+ 2\omega_3\lambda_{34}\lambda_{35}\lambda_{45} + \omega_3 \lambda_{35}^2 + \omega_4 \lambda_{45}^2 + \omega_5. 
\end{eqnarray*}
\end{example}

A classic result in the graphical models literature is the characterization of the conditional independence statements that hold for all densities belonging to the graphical model. This is determined by a graphical separation criterion, called \textit{$d$-separation}. 
\begin{definition}[Definition 2.3, \cite{trekSeparation2010}]
Let $A$, $B$ and $C$ be disjoint subsets of $[n]$. The set $C$ \textit{$d$-separates} $A$ and $B$ if every path (not necessarily directed) in $\cg$ connecting a node $i\in A$ to a node $j\in B$ contains a node $k$ that is either:
\begin{itemize}
\item a noncollider that belongs to $C$ or
\item a collider that does not belong to $C$ and has no descendants that belong to $C$,
\end{itemize}
where $k$ is a collider if there exist two consecutive edges $a \rightarrow k$ and $k\leftarrow b$ on the path.    
\end{definition}
The following result from \cite{LauritzenBook} ties the conditional independence statements of a DAG model to the $d$-separation statements of the DAG.

\begin{theorem}[Sec 3.2.2, \cite{LauritzenBook}]\label{theorem:d-separation theorem}
A set $C$ $d$-separates $A$ and $B$ in $\cg$ if and only if the conditional independence statement $X_A \indep X_B | X_C$ holds for every distribution in the graphical model associated to $\cg$.
\end{theorem}

Throughout this paper we will also often work with more general conditional independence models, which provide a direct generalization of graphical models. Let $\mathcal{C} = \{I_1 \indep J_1 | K_1, \ldots, I_\ell \indep J_\ell | K_\ell \}$ be a set of CI statements. Then the Gaussian \emph{conditional independence model} associated to $\mathcal{C}$ is the set of covariance matrices
\[
\cm_\mathcal{C} = \{\Sigma \in \pd_n ~:~ I \indep_\Sigma J | K \text{ for all } I \indep J | K \in \mathcal{C} \}. 
\]
These more general conditional independence models are naturally related to graphical models by the Hammersley-Clifford theorem \cite{algstat2018}, which states that the d-separation statements of $\cg$ yield the same distributions as those obtained from the structural equations of $\cg$. In other words, the conditional independence model defined by the global Markov property $\glob(\cg) = \{I \indep J | K ~:~ I \perp_\cg J | K \}$ is the same as the parameterized graphical model $\cm_\cg$ defined by the parameterization $\phi_\cg$. Note that throughout the paper we often follow the convention to write $I \indep J | K$ to mean that $X_I \indep X_J | X_K$.

Now, it is known that the normal random vector $X \sim \mathcal{N}(\mu, \Sigma)$ satisfies the conditional independence constraint $A \indep  B |C$ if and only if the submatrix $\Sigma_{A\cup C,B \cup C}$ has rank less than or equal to $|C|$.  We record this in an algebraic statement.

\begin{lemma}\label{lemma:minors}
    The CI statement $A \indep  B |C$ holds in a Gaussian random vector $X \sim \mathcal{N}(\mu, \Sigma)$ if and only if all minors of size $\#C+1$ of the submatrix $\Sigma_{A\cup C,B \cup C}$ are zero.
\end{lemma}

\begin{example}
Let $\cg$ be the DAG shown in \cref{figure:running-example}. The conditional independence statements that hold for this DAG are 
\[
1\indep 3 | \emptyset, \ 1 \indep 4 | \{2\}, \ 1\indep 5| \{2\}, \ \{1,2\} \indep 5 | \{3,4\}, \ 2\indep 3 | \emptyset.
\]
Now, if we consider the CI statement $\{1,2\} \indep 5 | \{3,4\}$, then by \cref{lemma:minors} we get that all the minors of size $3$ of the submatrix $\Sigma_{\{1,2,3,4\}\cup \{3,4,5\}}$ are zero. This implies that these minors are 
zero for every point in the model $\cm_\cg$. 
\end{example}


\section{The Algebraic Geometry of Unfaithful Distributions}
\label{sec: unfaithful}
In this section we explore the structure of the set of distributions in a graphical model that satisfy additional CI statements that are not implied by d-separation. In particular, we consider the subset of a Gaussian DAG model $\cm_\cg$ which satisfy an additional CI statement $i \indep j | K$ such that $i \not \perp_\cg j | K$ and show how tools from algebraic geometry can be used to study this problem. 


In \cite{TobiasThesis}, Boege outlines several key problems concerning conditional independence models one of which is the following well-studied problem. 

\begin{problem}[The CI Implication Problem]
\label{prob:CIimplication}
Let $\mathcal{C} = \{I_1 \indep J_1 | K_1, \ldots, I_\ell \indep J_\ell | K_\ell \}$ be a set of CI statements. For which disjoint sets $A, B, C \subseteq [n]$ does $\Sigma \in \mathcal{M}_\mathcal{C}$ imply $A \indep_\Sigma B | C$?\end{problem}

For arbitrary sets $\mathcal{C}$, this problem is extremely difficult and can only be solved in general with real quantifier elimination which can be doubly exponential in the number of variables involved \cite{TobiasThesis}. As previously mentioned, it has even been shown that no finite set of axioms can suffice to solve this problem in general \cite{studeny1992conditional, sullivant2009gaussian}. However, for certain sub-classes of conditional independence models this problem is easily solvable. For example, suppose that $\mathcal{C}$ consists of all the local Markov statements for a DAG $\cg$, meaning
\[
\mathcal{C} = \mathrm{local}(\cg) = \{ i \indep \nd(i)\setminus \pa(i) | \pa(i) ~:~ i \in V(\cg) \},
\]
where $\nd(i)$ is the set of \emph{non-descendants} of $i$ and $V(\cg)$ is the node set of $\cg$. Then it holds that $\cm_\mathcal{C} = \cm_\cg$ and thus one can easily check if an additional conditional independence statement $A \indep B | C$ holds for all $\Sigma \in \cm_\mathcal{C}$ by simply checking if $A \perp_\cg B | C$. The latter can be done in polynomial time \cite{koller2009probabilistic,pearl2009causality,spirtes2000causation}. One natural question is if there is a larger subclass of conditional independence models where \cref{prob:CIimplication} can also be solved efficiently. This leads us to the following question which will be the main focus of the remainder of this paper. 

\begin{problem}
\label{prob:CIimplicationOnDAGs}
Let $\cm_\cg$ be the graphical model on the DAG $\cg$ and suppose $i \not \perp_\cg j |K$. For which disjoint sets $A, B, C \subseteq [n]$ does $\Sigma \in \cm_\cg \cap \cm_{i \indep j | K}$ imply $A \indep_\Sigma B | C$? 
\end{problem}

In other words, we are studying Gaussian conditional independence models of the form $\mathcal{C} = \{I \indep J | K ~:~ I \perp_\cg J | K\} \cup \{i \indep j |K\}$ where $\cg$ is a DAG. In order to study this problem, we first introduce certain tools from algebraic geometry. The motivation behind adapting an algebraic perspective is that we can convert this problem into a principal ideal membership problem, which in turn becomes easy to solve in our setup.

Recall that the Gaussian DAG model is expressed as the image of the map $\phi_\cg$ as defined in \cref{definition:Gaussian DAG models}.
This parameterization can be rephrased in a combinatorial way with graphical sub-structures called \emph{treks}, as formalized by the so-called \emph{trek rule}; see, e.g., \cite{trekSeparation2010}.  In order to state the rule, we first give some basic definitions.

\begin{definition}
A \textit{trek} in $\cg$ from a node $i$ to a node $j$ is a pair $(P_L,P_R)$, where $P_L$ is a directed path from some node $s$ to $i$ and $P_R$ is a directed path from the same node $s$ to $j$. Here, $s$ is called the \textit{topmost} node of the trek, whereas $i$ and $j$ are called the \textit{leftmost} and \textit{rightmost} vertices of the trek, respectively. 

For any given trek $T=(P_L,P_R)$, the associated \textit{trek monomial} $m_T$ is given by 
\[
m_T= \lambda^L\omega_s \lambda^R=\omega_s \prod_{k\rightarrow l \in P_L}\lambda_{kl} \prod_{k\rightarrow l \in P_R} \lambda_{kl}.
\]
Note that in this setup, the edges in $P_L$ and $P_R$ are not necessarily disjoint, i.e., a variable $\lambda_{kl}$ can appear with multiplicity more than one in $m_T$. 
\end{definition}

With this terminology in place, we may formulate the \textit{trek rule} which states that the entries of  any covariance matrix $\Sigma \in \cm_\cg$ are given by
\[
\sigma_{ij}=\sum_T m_T= \sum_T \omega_{\text{topmost}(T)} \prod_{k\rightarrow l \in T}\lambda_{kl},
\]
where $T$ is the set of all possible treks between $i$ and $j$.  Now, let $\mathbb{R}[\Sigma] = \rr[\sigma_{ij} ~|~ 1 \leq i \leq j \leq n]$ and $\mathbb{R}[\lambda,\omega]$ be the rings of polynomials in the covariances $\sigma_{ij}$ and the non-zero entries of $(\Lambda,\Omega)$, respectively.  Then algebraically, the trek rule can be seen as a ring homomorphism $\phi^*_\cg$ from $\mathbb{R}[\Sigma] = \rr[\sigma_{ij} ~|~ 1 \leq i \leq j \leq n]$ to $\mathbb{R}[\lambda,\omega]$. 
Specifically,
\begin{eqnarray*}
\phi^*_\cg : \mathbb{R}[\Sigma] &\mapsto& \mathbb{R}[\lambda,\omega] \\
\sigma_{ij} &\mapsto& \sum_T m_T.
\end{eqnarray*}
Then $\phi^*_\cg(f)$ compactly denotes the result of substituting the trek rule expressions for each of the covariances appearing in the polynomial $f\in \mathbb{R}[\Sigma]$.

While the map $\phi_\cg$ naturally defines a statistical model, it can be extended to define an \emph{algebraic variety}. Associated to any subset $\cm\subseteq\pd_n$ is the ideal $I(\cm)$ that comprises all polynomials $f\in \mathbb{R}[\Sigma]$ that vanish on $\cm$.  In particular, we write $I_\cg:=I(\cm_\cg)$ for the  \textit{vanishing ideal} of the model $\cm_\cg$. In turn, associated to an ideal $I\subset\mathbb{R}[\Sigma]$, there is the algebraic variety $\cv(I) = \{\Sigma \in \rr^{\binom{n+1}{2}}~:~f(\Sigma) = 0 \text{ for all } f \in I\}$. This allows us to define the \emph{Zariski closure} of the graphical model $\cm_\cg$ which is $\overline{\cm_\cg} := \cv(I(\cm_\cg)).$ Throughout the next sections we will often work with the Zariski closure instead of the original model since we can apply algebraic tools to it naturally.

The model $\cm_\cg$ is parameterized and thus its Zariski closure forms an irreducible variety.  However, the same is not generally true for $\gmci = \cm_\cg \cap \cm_{i \indep j | K}$. Understanding the components of the algebraic variety $\vgmci:=\overline{\gmci}$ would typically require computing a \emph{primary decomposition} of the ideal $\mathcal{I}(\vgmci)$ which can be extremely difficult. However, the parameterization $\phi_\cg$ is injective and thus it is an isomorphism onto its image, $\cm_\cg$. This means that $\gmci$ is isomorphic to $\phi_\cg^{-1}(\gmci)$ and we can instead compute an irreducible decomposition of the Zariski closure of $\phi_\cg^{-1}(\gmci)$. This is advantageous since the vanishing ideal of  $\phi_\cg^{-1}(\gmci)$ is
\[
\phi_\cg^{\ast}(I_\cg + I_{i \indep j |K}) = 0 + \phi_\cg^{\ast}(|\Sigma_{iK, jK}|)
\]
which is a hypersurface so its primary decomposition can be computed via factorization. The following theorem summarizes this idea.

\begin{theorem}
\label{thm:DecomposeInEdgeSpace}
Let $\cg = (V, E)$ be a DAG and $i \indep j | K$ be a CI statement which is not implied by d-separation on $\cg$. Let $\phi_\cg^\ast$ be the ring map which corresponds to the parameterization of $\cm_\cg$ and $\phi_\cg^{\ast}(|\Sigma_{iK, jK}|) = \prod_{\ell = 1}^m f_\ell$ be an irreducible factorization. Then
\[
\gmci = \bigcup_{\ell = 1}^m \phi_\cg(\cv(f_\ell)) \cap \pd_n. 
\]    
\end{theorem}
\begin{proof}
Let $W_\ell = \overline{\phi_\cg(\cv(f_\ell))}$ be the Zariski closure of the image of $\cv(f_\ell)$ under $\phi_\cg$. First we show that each $W_\ell$ is an irreducible component of $\vgmci$. Note that $W_\ell \subset \vgmci$ by construction. Since $\cv(f_\ell)$ is an irreducible variety of codimension one in the parameter space $\cc^E \times \cc^V$, and $\phi_\cg$ is an injection, we immediately have that $W_\ell \subseteq \cc^{n(n+1)/2}$ is irreducible and $\dim(W_\ell) = \dim(V_\cg) - 1 = \dim(\vgmci)$. Thus each $W_\ell$ must be an irreducible component of $\vgmci$ since it is an irreducible subvariety with the highest possible dimension. 

Now suppose that $\vgmci = \cup_\alpha V_\alpha$ is a decomposition into irreducible varieties. Thus $\gmci = \vgmci \cap \pd_n = \cup_{\alpha} V_\alpha \cap \pd_n$. So to complete the proof it suffices to show that if $V_\alpha \cap \pd_n \neq \emptyset$, then $V_\alpha \cap \pd_n \subseteq W_\ell$ for some $\ell$. Since $\phi_\cg$ is injective, it is an isomorphism when restricted to its image. Thus we have that $\phi_\cg^{-1}(V_\alpha \cap \pd_n) \subseteq V(\phi_\cg^\ast(|\Sigma_{iK, jK}|))$ by construction but $\overline{\phi_\cg^{-1}(V_\alpha \cap \pd_n)}$ must be irreducible so there exists some $\ell \in [m]$ such that $\phi_\cg^{-1}(V_\alpha \cap \pd_n) \subset \cv(f_\ell)$ which immediately implies that $V_\alpha \cap \pd_n \subseteq W_\ell$.
\end{proof}

The above theorem makes it significantly easier to characterize the components of $\gmci$ (by which we mean the components of $\vgmci$ which have non-empty intersection with $\pd_n$) since they correspond to the factors of the polynomial $\phi_\cg^{\ast}(|\Sigma_{iK, jK}|)$. This is demonstrated by the following example.

\begin{example}
\label{ex:DecomposeIntoGM}
Let $\cg$ be the graph pictured in \cref{figure:running-example}. Observe that $1 \not \perp 5 | 4$ in this graph since the path $1 \to 2 \to 4 \leftarrow 3 \to 5$ is d-connecting with respect to the conditioning set $K = \{4\}$. So in this case our model of interest is $\cm_{\cg, 1 \indep 5 |4} = \cm_\cg \cap \cm_{1 \indep 5 | 4}$. By \cref{thm:DecomposeInEdgeSpace}, we have that the irreducible components of $\gmci$ correspond to the factors of $\phi_\cg^\ast(|\Sigma_{14, 45}|)$. Using the \emph{trek rule} described in \cref{sec:SEMBackground}, $\phi_\cg(\sigma_{ij})$ can be easily computed and thus so can $\phi_\cg^\ast(|\Sigma_{14, 45}|)$. Doing so yields
\[
\phi_\cg^\ast(|\Sigma_{14, 45}|) = \phi_\cg^\ast(\sigma_{15}\sigma_{44} - \sigma_{14}\sigma_{45})= \omega_1 \omega_3 \lambda_{12} \lambda_{24} \lambda_{34} \lambda_{35}. 
\]
Observe that $\phi_\cg(\cv(\omega_i)) \cap \pd_n = \emptyset$ since for any matrix $\Sigma \in \phi_\cg(\cv(\omega_i))$, it holds that
\[
|\Sigma| = |(I-\Lambda)^{-T} \Omega (I - \Lambda)^{-1}| = |\Omega| = 0. 
\]
On the other hand, the remaining four factors of $\phi_\cg^\ast(|\Sigma_{14, 45}|)$ are simply the variables $\lambda_{12}, \lambda_{24}, \lambda_{34}, \lambda_{35}$. 
Recall that $\cv(\lambda_{ij}) = \{\Lambda \in \rr^{E} ~|~ \lambda_{ij} = 0\}$. Thus if $i \to j$ is an edge in $\cg$, then $\phi_\cg(\cv(\lambda_{ij})) \cap \pd_n = \cm_{\cg \setminus i \to j}$ where $\cg \setminus i \to j$ is the DAG obtained by deleting the edge $i \to j$ from $\cg$. In other words, the components of the form $\phi_\cg(V(\lambda_{ij})) \cap \pd_n$ are exactly the graphical sub-models of $\cm_\cg$ which correspond to deleting the edge $i \to j$ from the original graph $\cg$.  
So we see in this case that 
\[
\cm_{\cg, 1 \indep 5 |4} \;=\; \cm_{\cg \setminus 1 \to 2} \cup \cm_{\cg \setminus 2 \to 4} \cup \cm_{\cg \setminus 3 \to 4} \cup \cm_{\cg \setminus 3 \to 5}.
\]
\end{example}

One immediate consequence of \cref{thm:DecomposeInEdgeSpace} is that solving the conditional independence implication problem for models of the form $\gmci$ is actually quite tractable from an algebraic perspective. The following corollary gives a sufficient condition for which other CI statement hold for all distributions in the model $\gmci$.  It is formulated using the concept of \emph{saturation} \cite{coxlittleoshea}.

\begin{corollary}
\label{cor:CIimplicationViaPrincipal}
Let $a, b \in [n]$ be disjoint elements and $C \subseteq [n] \setminus \{a, b\}$. Then $a \indep_\Sigma b | C$ for all distributions $\Sigma \in \gmci$ if
\[
\phi_\cg^{\ast}(|\Sigma_{aC, bC}|) \;\in\; \langle \phi_\cg^\ast(|\Sigma_{iK, jK}|) \rangle : \prod_{S \subseteq V(\cg)} \phi_\cg^\ast(|\Sigma_{S, S}|)^\infty.
\]
\end{corollary}

\begin{example}
Let $\cg$ be the graph pictured in \cref{figure:running-example}, and consider the additional statement $i \indep j | K = 1 \indep 2 | 5$. Then
\[
\langle \phi_\cg^\ast(|\Sigma_{iK, jK}|\rangle =\langle  \phi_\cg^\ast(|\Sigma_{15, 25}|) \rangle = \langle \lambda_{24}(\omega_1 \lambda_{12}^2 + \omega_2)(\omega_3 \lambda_{34} \lambda_{35} \lambda_{45} + \omega_3 \lambda_{35}^2 + \omega_5) \rangle . 
\]
The saturation described in the previous lemma is with respect to the principal minors of $\Sigma$, which simply means that we remove any factor of the above polynomial that corresponds to a principal minor. In this case, the second factor corresponds to the principal minor $\Sigma_{2, 2}$ while the other two factors do not correspond to any principal minor, which means
\[
\langle \phi_\cg^\ast(|\Sigma_{iK, jK}|) \rangle : \prod_{S \subseteq V(\cg)} \phi_\cg^\ast(|\Sigma_{S, S}|)^\infty = \langle \lambda_{24}(\omega_3 \lambda_{34} \lambda_{35} \lambda_{45} + \omega_3 \lambda_{35}^2 + \omega_5) \rangle. 
\]
Now suppose that we want to check if $1 \indep 4 |5$ for all covariance matrices in $\cm_{\cg, 1\indep 2 | 5}$. Then by the previous corollary we simply need to compute
\[
\phi_\cg^\ast(|\Sigma_{15, 45}|) = \omega_1 \lambda_{12} \lambda_{24} (\omega_3 \lambda_{34} \lambda_{35} \lambda_{45} + \omega_3 \lambda_{35}^2 + \omega_5).
\]
It is clear that in this case $\phi_\cg^\ast(|\Sigma_{15, 25}|)$ divides $\phi_\cg^\ast(|\Sigma_{15, 45}|)$ thus it holds that $1 \indep 4 |5$ for all $\Sigma \in \cm_{\cg, 1\indep 2 | 5}$. 
\end{example}

The above corollary provides a relatively easy algebraic condition for testing CI implication. This is because the ideal  $\langle \phi_\cg^\ast(|\Sigma_{iK, jK}|) \rangle : \prod_{j} \phi_\cg^\ast(|\Sigma_{\pa(j), \pa(j)}|)^\infty$ is a \emph{principal} ideal obtained by removing any factor of $\phi_\cg^\ast(|\Sigma_{iK, jK}|)$ which corresponds to one of the principal minors $\phi_\cg^\ast(|\Sigma_{\pa(j), \pa(j)}|)$; more details on this saturation operation can be found in \cite[Section 4]{algstat2018}. Thus to test if $\phi_\cg^{\ast}(|\Sigma_{aC, bC}|)$ belongs to it one simply needs to test if the single generator of  $\langle \phi_\cg^\ast(|\Sigma_{iK, jK}|) \rangle : \prod_{j} \phi_\cg^\ast(|\Sigma_{\pa(j), \pa(j)}|)^\infty$ divides $\phi_\cg^\ast(|\Sigma_{aC, bC}|)$ which can be done in polynomial time in the number of variables. However, the corollary does have two notable drawbacks. First, it does still require that one compute $\phi_\cg^{\ast}(|\Sigma_{aC, bC}|)$ and $\phi_\cg^{\ast}(|\Sigma_{iK, jK}|)$ which becomes much more difficult as the size of $\cg$ grows. Second, it provides a sufficient condition for testing CI implication but it is not necessary. This is because \cref{cor:CIimplicationViaPrincipal} is actually testing which CI statements $a \indep_\Sigma b | C$ are true  for all matrices $\Sigma \in \vgmci$. However, it can be that there are CI statements which hold for all positive definite matrices $\Sigma \in \gmci$ but not all matrices in the algebraic closure $\vgmci$ as shown in \cite[Example 3.3]{TobiasThesis}.

\section{Decomposition Into Graphical Models}\label{sec:graphicalImplication}
While solving \cref{prob:CIimplication} is extremely difficult, we've seen that \cref{prob:CIimplicationOnDAGs} is much more tractable though it still might require computing determinants over a polynomial ring which becomes infeasible for large matrices. As we have already discussed, \cref{prob:CIimplication} is significantly easier for graphical models since a given conditional independence holds if and only if the corresponding d-separation holds in the graph which can be easily checked. While the model $\gmci$ is no longer a graphical model, it is possible that it can be written as a union of graphical models as we saw in \cref{ex:DecomposeIntoGM}. We are particularly interested in the case where $\gmci$ decomposes into a union of graphical models since the conditional independence implication problem can also be easily solved with d-separation in this setting. The following proposition makes this explicit.  

\begin{proposition}
\label{prop:ImplicationForGMUnion}
Let $\cg = (V, E)$ be a DAG and $i \indep j | K$ be a CI statement which is not implied by d-separation on $\cg$. If $\gmci = \cup_{\ell = 1}^m \cm_{\cg_\ell}$ for some DAGs $\cg_1, \ldots, \cg_m$ then $A \indep_\Sigma B | C$ for all $\Sigma \in \gmci$ if and only if $A \perp_{G_\ell} B |C$ for all $\ell \in [m]$. 
\end{proposition}
\begin{proof}
Observe that $A \indep_\Sigma B |C$ for all $\Sigma \in \gmci$ if and only if $A \indep_\Sigma  B | C \in \cm_{\cg_\ell}$ for all $\ell \in [m]$. Since $\cm_{\cg_\ell}$ is a Gaussian linear SEM on $\cg$, $A \indep_\Sigma B | C$ if and only if $A \perp_{\cg_\ell} B | C$ by \cref{theorem:d-separation theorem}, which completes the proof. 
\end{proof}

We can check whether $A \perp_{\cg_\ell} B | C$ holds in polynomial time in the number of vertices $n$ \cite{van2019separators} and therefore \cref{prop:ImplicationForGMUnion} implies that we can solve \cref{prob:CIimplication} in polynomial time when $\gmci = \cup_{\ell = 1}^m \cm_{\cg_\ell}$ by simply checking d-separations in each graph $\cg_\ell$. The next proposition gives a necessary and sufficient algebraic condition for when $\gmci$ decomposes into a union of graphical models. 

\begin{proposition}
\label{proposition: monomial}
Let $\cg = (V, E)$ be a DAG and $i \indep j | K$ be a CI statement which is not implied by d-separation on $\cg$. Then there exist graphs $\cg_1, \cg_2, \ldots, \cg_m$ such that
\[
\gmci = \cup_{\ell = 1}^m \cm_{\cg_\ell}
\]
if and only if $\phi_\cg^\ast(|\Sigma_{iK, jK}|)$ is a monomial in the ring $\rr[\Lambda, \Omega]$. 
\end{proposition}
\begin{proof}
Let $\gmci = \bigcup_{\ell = 1}^m \phi_\cg(\cv(f_\ell)) \cap \pd_n$ be a decomposition into irreducible components as guaranteed by \cref{thm:DecomposeInEdgeSpace}. Observe that if $\phi_\cg^\ast(|\Sigma_{iK, jK}|) = \prod_{\ell = 1}^m f_\ell$ is a monomial, then each $f_\ell = \lambda_{ij}$ for some edge $i \to j \in \cg$ or $f_\ell = \omega_i$. If $f_\ell = \lambda_{ij}$ then $\phi_\cg(\cv(f_\ell)) \cap \pd_n = \cm_{\cg \setminus i \to j}$. On the other hand, if $f_\ell = \omega_i$ then $\phi_\cg(\cv(f_\ell)) \cap \pd_n = \emptyset$.

Now suppose that $\phi_\cg^\ast(|\Sigma_{iK, jK}|) = \prod_{\ell = 1}^m f_\ell$ is not a monomial so there exists some $\ell$ such that $f_\ell$ is not a monomial. Now observe that $\cm_{\cg'} \subseteq \cm_{\cg}$ if and only if $\cg'$ is Markov equivalent to a subgraph of $\cg$. Moreover, if $\dim(\cm_{\cg'}) = \dim(\cm_\cg) - 1$, then $\cm_{\cg'} = \cm_{\cg \setminus (i, j)}$ for some edge $(i, j) \in \cg$. However, this means that $\overline{\phi_\cg^{-1}(\cm_{\cg'})} = \cv(\lambda_{ij}) \neq V(f_\ell)$ by assumption.   
\end{proof}

While the previous proposition provides an algebraic condition for when $\gmci$ decomposes into graphical models, it can be difficult to check in practice since it involves computing $\phi_\cg^\ast(|\Sigma_{iK, jK}|)$ which becomes expensive as the size of $\cg$ grows. We now focus on developing an equivalent graphical condition which is easy to check. In order to do this, we first analyse the structure of $\phi_\cg(|\Sigma_{A,B}|)$, where $A$ and $B$ are arbitrary sets of the same size. This scenario was studied in \cite{positivity2013}, where the authors developed the concept of \textit{treks systems} and \textit{no sided intersection}. 

\begin{definition}
Let $A$ and $B$ be sets of $k$ vertices.
\begin{enumerate}

\item A \textit{trek system} $\bold{T}$ from $A$ to $B$ consists of $k$ treks whose leftmost vertices exhaust the set $A$ and whose rightmost vertices exhaust the set $B$. The \textit{trek system monomial} $m_{\bold{T}}$ is defined as the the product of the trek monomials $m_T$ ranging over $T \in \bold{T}$.

\item A trek system $\bold{T}$ from $A$ to $B$ is said to have \textit{no sided intersection} if the left parts $P_L$ for $(P_L,P_R) \in \bold{T}$ are mutually node-disjoint and also the right parts $P_R$ for $(P_L, P_R) \in \bold{T}$ are mutually node-disjoint (although any $P_L$ may have vertices in common with any $P_{R'}$). We use $\mathcal{T}(A, B)$ to denote set of trek systems from $A$ to $B$ with no sided intersection.
\end{enumerate}
\end{definition}
    
Any trek system $\mathbf{T}$ provides a bijection between the sets $A$ and $B$. For a given linear ordering of $A$ and $B$, this bijection can be determined by a permutation $\pi$ on $n$ elements. Thus the sign of the trek system $\mathbf{T}$ is defined as  $\mathrm{sign}(\mathbf{T}) := \mathrm{sign}(\pi)$. We now state the relevant result from \cite{positivity2013} which we will use to derive our new graphical condition. 

\begin{corollary}\label{corollary:det equivalence}[Corollary 3.5, \cite{positivity2013}]
Let $G = (V, D)$ be a DAG and $A$ and $B$ be subsets of $V$ of the same cardinality. Then 
\[
|\Sigma_{A,B}|= \sum_{[\bold{T}]_{\sim} \in \mathcal{T}(A,B)/ \sim}   \text{sign}(\bold{T})2^{|UD(T)|}m_{\bold{T}}
\]
where the sum runs over equivalence classes of the relation $\sim$ defined by $\bold{T} \sim \bold{T'}$ if and only if $m_{\bold{T}} = m_{\bold{T'}}$.
\end{corollary}

Here $UD(T)$ corresponds to the number of \textit{up-down} cycles in $\bold{T}$. The concept of up-down cycle is explained in \cref{example:updowncycles}. In general, an up-down cycle consists of two source nodes $a,b$, two sink nodes $a_1,b_1$ and two pairs of directed paths from $a$ and $b$ to $a_1$ and $b_1$, respectively. The idea behind this construction is that we can get two different trek systems with no sided intersection by collectively using the same set of edges. More details about the construction can be found in \cite{positivity2013}.

\cref{corollary:det equivalence} immediately gives us a necessary and sufficient condition on when $\phi(|\Sigma_{iK,jK}|)$ is a monomial, i.e., either when there exists a unique trek system between $\{i,K\}$ and $\{j,K\}$ with no sided intersection, or if every trek system with no sided intersection has the same trek system monomial. In the next Corollary, we show that the latter is not possible in our setup.

\begin{corollary}\label{corollary:monomial condition}
Let $G=(V,E)$ be a DAG and $\{i\},\{j\},K \subset V$ be set of vertices such that $K$ does not $d$-separate $i$ from $j$. Then $\phi( |\Sigma_{iK,jK}|)$ is a monomial if and only if there is a unique trek system $\bold{T}$ from $\{i,K\}$ to $\{j,K\}$ with no sided intersection.    
\end{corollary}
\begin{proof}
It is clear from \cref{corollary:det equivalence} that if there exists a unique trek system from $\{i,K\}$ to $\{j,K\}$ with no sided intersection, then the image of $|\Sigma_{\{i,K\},\{j,K\}}|$ is a monomial (as the sum runs over a single equivalence class of size one). Thus, in order to prove the only-if direction, we need to show that if there exist multiple trek systems with no sided intersection having the same trek system monomial, then there must exist another trek system having a different trek system monomial, making the image of $|\Sigma_{\{i,K\},\{j,K\}}|$ a non monomial.

So, let $\bold{T}$ and $\bold{T'}$ be two different trek systems with no sided intersection with $m_{\bold{T}}=m_{\bold{T'}}$. We have two possible scenarios in this case: either $\bold{T}$ (and $\bold{T'}$) contain a trek between $i$ and $j$, or the trek system(s) is of the form  $i \leftarrow \ldots \rightarrow k_1, k_1 \leftarrow \ldots \rightarrow k_2, \ldots, k_s \leftarrow \ldots \rightarrow j, k_{a_1}\leftrightarrow k_{b_1}, \ldots, k_{a_n}\leftrightarrow k_{b_n}$. In the first case, we can construct new trek system $\bold{T''}$ having a different trek system monomial by taking the trek between $i$ and $j$ along with the empty treks at $k_l$. In the second case, we can assume without loss of generality that the up-down cycle is contained in the first two treks. This implies that the first two treks $\bold{T}$ are of the following form:
\begin{eqnarray*}
i\leftarrow \ldots \leftarrow x_0 \leftarrow \ldots  \leftarrow s_1 \rightarrow \ldots x_1 \rightarrow \ldots \rightarrow k_1, \\
k_1\leftarrow \ldots \leftarrow x_1 \leftarrow \ldots \leftarrow s_2 \rightarrow \ldots x_0 \rightarrow \ldots \rightarrow k_2.
\end{eqnarray*}
The up-down cycle in this case is $x_0 \leftarrow \ldots \leftarrow s_1 \rightarrow \ldots \rightarrow x_1 \leftarrow \ldots \leftarrow s_2 \rightarrow \ldots \rightarrow x_0$. Using this structure, we create the new trek system $\bold{T''}$ by replacing only the first two treks in $\bold{T}$ with the following two treks:
\begin{eqnarray*}
 i\leftarrow \ldots  \leftarrow x_0  \rightarrow \ldots \rightarrow k_2, k_1 \leftrightarrow k_1. 
\end{eqnarray*}
In other words, we create $\bold{T''}$ by omitting the up-down cycle from $\bold{T}$, which also implies that $m_{\bold{T}}$ is not equal to $m_{\bold{T''}}$.
\end{proof}

\begin{figure}
\centering
\begin{tikzpicture}[>=stealth',shorten >=1pt,auto,node distance=0.8cm,scale=.9, transform shape,align=center,minimum size=3em]
\node[state] (w5) at (0,0) {$5$};
\node[state] (w7) at (3,0) {$7$};
\node[state] (w8) at (6,0) {$8$};
\node[state] (w9) at (9,0) {$9$}; 
\node[state] (w3) at (1.5,2.5) {$3$};
\node[state] (w6) at (4.5,1.5) {$6$};
\node[state] (w4) at (7.5,2.5) {$4$};
\node[state] (w1) at (1.5,4.5) {$1$};
\node[state] (w2) at (7.5,4.5) {$2$};


\path[->]   (w1) edge    (w3);
\path[->]   (w1) edge    (w4);
\path[<-]   (w3) edge    (w2);
\path[->]   (w2) edge    (w4);
\path[->]   (w3) edge    (w5);
\path[->]   (w3) edge    (w7);
\path[->]   (w6) edge    (w7);
\path[->]   (w6) edge    (w8);
\path[->]   (w4) edge    (w8);
\path[->]   (w4) edge    (w9);

\end{tikzpicture}
\caption{A DAG $\cg$ with an up-down cycle.}
\label{figure:up-down cycle}
\end{figure}
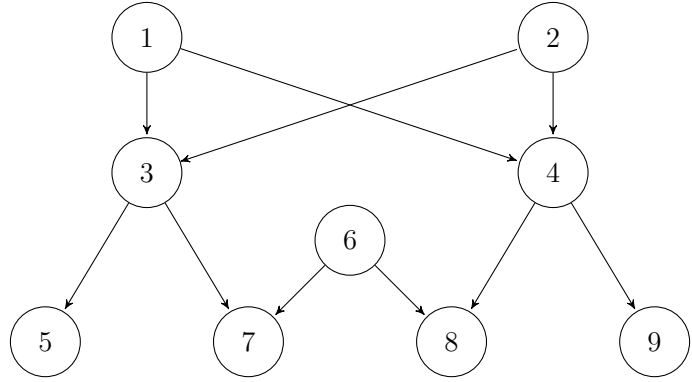

\begin{example}\label{example:updowncycles}
Consider the DAG in \cref{figure:up-down cycle}. Let $A$ and $B$ be the sets $\{5,7,8\}$ and $\{7,8,9\}$, respectively. Any trek system between $A$ and $B$ is a collection of three treks whose leftmost vertices exhaust $\{5,7,8\}$ and the rightmost vertices exhaust $\{7,8,9\}$. The following is a trek system between $A$ and $B$ with no sided intersection:
\[
\bold{T}=\{5\leftarrow 3 \rightarrow 7, 7\leftarrow 6 \rightarrow 8, 8\leftarrow 4 \rightarrow 9\}.
\]
Observe that there is no intersection between the left parts (and similarly right parts) of the treks. The corresponding trek system monomial is $m_{\bold{T}}=\omega_3\omega_4\omega_6\lambda_{35}\lambda_{37}\lambda_{67}\lambda_{68}\lambda_{48}\lambda_{49}$. 

The DAG $\cg$ has an up-down cycle $1\rightarrow 3 \leftarrow 2 \rightarrow 4 \leftarrow 1$. Using this cycle, we construct two more trek systems with no sided intersection having the trek system monomial. Let $\bold{T_1}$ and $\bold{T_2}$ be the following two trek systems:
\begin{eqnarray*}
   \bold{T_1}&=& \{5\leftarrow 3 \leftarrow 1 \rightarrow 4 \rightarrow 9, \hspace{.2cm} 8\leftarrow 4 \leftarrow 2 \rightarrow 3 \rightarrow 7, \hspace{.2cm} 7\leftarrow 6 \rightarrow 8\}, \\
   \bold{T_2}&=& \{5\leftarrow 3 \leftarrow 2 \rightarrow 4 \rightarrow 9, \hspace{.2cm} 8\leftarrow 4 \leftarrow 1 \rightarrow 3 \rightarrow 7, \hspace{.2cm} 7\leftarrow 6 \rightarrow 8\}.
\end{eqnarray*}
One can check that both trek systems have no sided intersection and also have the same trek system monomial as they collectively have the same topmost vertices and same edges (with same multiplicity as well).
\end{example}

While the previous corollary provides a necessary and sufficient condition for when $\gmci$ decomposes into graphical models, we would ideally like a condition which is more directly checkable given the graph $\cg$ and the CI statement $i \indep j | K$. Recall that the vertices $i$ and $j$ are said to be \textit{$d$-connected} given the set $K$ if there exists a trek between $i$ and $j$ which does not contain any node from $K$, or if there exists a path where every collider lies in $K$. In other words, $i$ and $j$ are said to be $d$-connected given $K$ if they are not $d$-separated given $K$. Now, the trek systems between $\{i,K\}$ and $\{j,K\}$ with no sided intersection can also be seen as $d$-connecting paths between $i$ and $j$ given $K$. We prove this correspondence in the following lemma. 

\begin{lemma}\label{lemma:trek system-d connection}
Let $\mathcal{T}(\{i,K\},\{j,K\})$ be the set of trek systems between $\{i,K\}$ and $\{j,K\}$ with no sided intersection and $S(i,j)_K$ be the set of $d$-connecting paths between $i$ and $j$ given $K$. Then each trek system in $\mathcal{T}(\{i,K\},\{j,K\})$ corresponds to a unique $d$-connecting path in $S(i,j)_K$. Further, each $d$-connecting path in $S(i,j)_K$ also corresponds to at least one trek system in $\mathcal{T}(\{i,K\},\{j,K\})$. 
\end{lemma}
\begin{proof}
Let $i<j$ and $K=\{k_1,k_2,\ldots,k_m\}$. We first show that any trek system in $\mathcal{T}(\{i,K\},\{j,K\})$ has the following form:
\begin{enumerate}
    \item $i\leftarrow \ldots \rightarrow j, k_{a_1}\leftrightarrow k_{b_1}, k_{a_2}\leftrightarrow k_{b_2}, \ldots, k_{a_n}\leftrightarrow k_{b_n}$, where $1\leq a_l,b_l \leq m$ and the trek between $i$ and $j$ doesn't contain any node in $K$,
    \item $i \leftarrow \ldots \rightarrow k_1, k_1 \leftarrow \ldots \rightarrow k_2, \ldots, k_s \leftarrow \ldots \rightarrow j, 
    k_{a_1}\leftrightarrow k_{b_1}, \ldots, k_{a_n}\leftrightarrow k_{b_n}$ for all $a_l,b_l \geq s+1$.
\end{enumerate}

In both cases, we can have $a_l$ and $b_l$ as equal. In such situations, we can either consider the empty trek $k_{a_l}\leftrightarrow k_{a_l}$, or take any non simple trek of the form $k_{a_l} \leftarrow p \rightarrow k_{a_l}$ with $p<k_{a_l}$ in the trek system. Although this structure will be useful in the next Proposition, the proof of the current lemma is independent of this structure.

To show this, we look at all the treks with $i$ as the leftmost node in the trek systems. If any trek system contains a trek between $i$ and $j$ (i.e., with $j$ as its rightmost node), then this trek cannot contain any node $k_l\in K$ as that would form a sided intersection. Thus, the system has to be of type $1$. Similarly, if $i \leftarrow \ldots \rightarrow k_1$ is a trek in one of the trek systems, then there must exist a trek with $k_1$ as the leftmost node in order to exhaust the set $\{i,K\}$. However, this trek cannot have $k_1$ as its topmost node as that would cause a sided intersection with $i \leftarrow \ldots \rightarrow k_1$. Thus, the second trek has to be of the form $k_1 \leftarrow \ldots \rightarrow k_2$. Iterating the same argument, we can conclude that the system has to be of type $2$.

Now that we have shown the types of trek systems in $\mathcal{T}(\{i,K\},\{j,K\})$, it is easy to see the corresponding $d$-connecting paths. For trek systems of type $1$, the corresponding $d$-connecting path is the trek $i\leftarrow \ldots \rightarrow j$ as it does not contain any $k_l\in K$. Similarly, for trek systems of type $2$, the corresponding $d$-connecting path is obtained by taking the union of the treks $i \leftarrow \ldots \rightarrow k_1 \cup  k_1 \leftarrow \ldots \rightarrow k_2 \cup \ldots \cup k_s \leftarrow \ldots \rightarrow j$. As all the colliders in this path lie in $K$, it is indeed $d$-connecting given $K$.

The reverse correspondence is trivial as any path in $S(i,j)_K$ is either a trek between $i$ and $j$ not containing any $k_l$ (which corresponds to the trek systems of type $1$), or is a path which has colliders at each $k_l$ that lies in it (corresponding to trek systems of type $2$).
\end{proof}

Using the above correspondence, \cref{corollary:monomial condition} can also be stated in terms of $d$-connecting paths in the following way:
\begin{theorem}\label{theorem:d-connection monomial}
Let $G=(V,E)$ be a DAG and $\{i\},\{j\},K \subset V$ be sets of vertices such that $K$ does not $d$-separate $i$ from $j$. Then $\phi(\det \Sigma_{iK,jK})$ is a monomial if and only if
\begin{itemize}
\item there is a unique $d$-connecting path from $i$ to $j$ given $K$ (denoted by $i \leftrightarrow j$), and
\item For all $k_\ell \notin i \leftrightarrow j$, it holds that $\pa(k_\ell) \subseteq i\leftrightarrow j \cup K$
\end{itemize}
\end{theorem}
\begin{proof}
As each $d$-connecting path between $i$ and $j$ given $K$ corresponds to at least one trek system between $\{i,K\}$ and $\{j,K\}$ with no sided intersection (as seen in \cref{lemma:trek system-d connection}), it is necessary to have a unique $d$-connecting path in order to have a monomial image of $|\Sigma_{iK, jK}|$. However, this condition is not sufficient. This is because, for each $d$-connecting path the trivial trek system that corresponds to it is the one that has empty treks for each $k_l \in K$ that do not lie in the path. However, if any $k_l \in K$ that does not lie in the path has an incoming edge of the form $p\rightarrow k_l$ (with $p\in V\setminus i\leftrightarrow j$), then we can use the non simple trek $k_l \leftarrow p \rightarrow k_l$ to construct another trek system. (Note that this non simple trek would form a sided intersection with $i \leftrightarrow j$ if $p\in i\leftrightarrow j$.) Thus, requiring the second condition forces $k_l$ to be the topmost node for any trek in $V\setminus i\leftrightarrow j$ that contains it, implying that there cannot exist another trek system without any sided intersection apart from the one containing empty treks for $k_l$.     
\end{proof}

\cref{corollary:monomial condition} gives us a sufficient graphical condition for when $\gmci$ corresponds to a graphical model. As we saw earlier, each of these components correspond to some $\lambda_{ij}$ with $(i,j)\in E$, and so the graphical models which appear in the decomposition are precisely the ones obtained after deleting the edge $(i,j)$ from $\cg$. Now, when the condition in \cref{corollary:monomial condition} is not satisfied, i.e., when there are multiple $d$-connecting paths between $i$ and $j$ given $K$, the components obtained after primary decomposition are not necessarily graphical models. However, some of the components of the decomposition may correspond to graphical models in certain cases. We provide a graphical condition for such scenarios in the corollary below. 

\begin{corollary}\label{corollary:graphical component}
Let $\cg=(V,E)$ be a DAG and $\{i\},\{j\},K \subset V$ be set of vertices such that $K$ does not $d$-separate $i$ from $j$. If every $d$-connecting path between $i$ and $j$ given $K$ contains the edge $(a,b)$, then there exists an irreducible component of $\gmci$ which corresponds to the graphical model $\cg\setminus a\rightarrow b$.   
\end{corollary}
\begin{proof}
As every $d$-connecting path between $i$ and $j$ given $K$ passes through the edge $(a,b)$, the variable $\lambda_{ab}$ appears in every trek system monomial from $\{i,K\}$ to $\{j,K\}$ with no sided intersection. Thus, by \cref{corollary:det equivalence} we can conclude $\phi (|\Sigma_{iK,jK}|)$ factorizes into $\lambda_{ab}\prod_{l=1}^m f_l$. Thus by \cref{thm:DecomposeInEdgeSpace}, one irreducible component of $\gmci$ is $\phi_\cg(V(\lambda_{ab}))$ which corresponds to the graphical model $\cg \setminus a\rightarrow b$. 
\end{proof}

The following example illustrates the previous lemmas. 

\begin{example}\label{example:graphical component}
Consider the DAG in \cref{figure:running-example}. The CI statements that hold for this DAG are 
\[
1\indep 3,~ 2\indep 3,~ 1 \indep 4 |2,~ 1\indep 5| 2,~ 1\indep 5 | \{3,4\},~  2\indep 5 | \{3,4\}.
\]
We again consider the model $\cm_{\cg, 1 \indep 5 |4} = \cm_\cg \cap \cm_{1 \indep 5 |4}$. Observe that there is exactly one $d$-connecting path between $1$ and $5$ given $4$, which is the path $1 \rightarrow 2 \rightarrow 4 \leftarrow 3 \rightarrow 5$. As there is no other node in the conditioning set apart from $4$, the conditions in \cref{theorem:d-connection monomial} are satisfied. Thus, we know that $\phi(|\Sigma_{\{1,4\},\{4,5\}}| )$ is a monomial. As we saw in \cref{ex:DecomposeIntoGM}, computing the image yields
\[
\phi_\cg^\ast(|\Sigma_{\{1,4\},\{4,5\}}| )=\phi(\sigma_{14}\sigma_{45}-\sigma_{15}\sigma_{44})=\omega_1 \omega_3 \lambda_{12}\lambda_{24}\lambda_{34}\lambda_{35},
\]
which is precisely the trek monomial corresponding to the $d$-connecting path. 

Now, consider adding the statement $1\indep 2 |5$ to $\cm_\cg$. Observe that even though there exists a unique $d$-connecting path between $1$ and $2$ given $5$ (which is just the edge $1\rightarrow 2$), the node in the conditioning set (i.e., $5$) has incoming edges of the form $3\rightarrow 5$ and $4\rightarrow 5$. Thus, by \cref{theorem:d-connection monomial}, we know that $\phi(|\Sigma_{\{1,5\},\{2,5\}}| )$ is not a monomial. Computing the image gives us that
\begin{eqnarray*}
\phi_\cg^\ast(| \Sigma_{\{1,5\},\{2,5\} }|) = \phi(\sigma_{12}\sigma_{55}-\sigma_{15}\sigma_{25})&=&\omega_1 \omega_3 \lambda_{12}\lambda_{34}^2\lambda_{45}^2 + 
2\omega_1 \omega_3 \lambda_{12}\lambda_{34}\lambda_{35}\lambda_{45}\\
&&+ 
\omega_1 \omega_3 \lambda_{12}\lambda_{35}^2+
\omega_1 \omega_4 \lambda_{12}\lambda_{45}^2 +
\omega_1 \omega_5\lambda_{12}.
\end{eqnarray*}
However, every $d$-connecting path between $1$ and $2$ given $5$ does contain the edge $1\rightarrow 2$ (which is trivial in this case). Thus, by \cref{corollary:graphical component}, we know that $\cm_{\mathcal{G}\setminus 1\rightarrow 2}$ must be an irreducible component of $\cm_{\cg, 1\indep2
|5}$. This is indeed the case, as we can see from the computation that $\lambda_{12}$ is an irreducible factor of $\phi_\cg^\ast(|\Sigma_{\{1,5\},\{2,5\} }|)$. 
\end{example}

\begin{remark}
The technique developed in this section can be implemented iteratively to potentially determine if  $\mathcal{M}_\cg \cap \cm_\mathcal{C}$ decomposes into a union of graphical models for any arbitrary collection of additional CI statements $\mathcal{C}$. For instance, if we add the statements $1 \indep 5|\{4\}$ and $1\indep 2 | \{5\}$ to $\cm_\cg$ as seen in \cref{example:graphical component}, then the irreducible components of the model can be determined by applying the above results twice. Thus, in order to obtain an irreducible decomposition of 
$\cm_\cg \cap \cm_{1 \indep 5 | 4} \cap \cm_{1 \indep 2 |5}$
we first analyze the components of $\cm_\cg  \cap \cm_{1 \indep 5 |4}$. Recall that this model decomposed into
\[
\cm_{\cg, 1 \indep 5 |4} = \cm_{\cg \setminus 1 \to 2} \cup \cm_{\cg \setminus 2 \to 4} \cup \cm_{\cg \setminus 3 \to 4} \cup \cm_{\cg \setminus 3 \to 5}.
\]
Now, observe that the component $\cm_{\cg\setminus 1\rightarrow 2}$ already satisfies the statement $1\indep 2 | \{5\}$, whereas the other three components do not; however, the other three models all satisfy 
\[
\cm_{\cg\setminus i \to j} \cap \cm_{1 \indep 2 | 5} \subseteq \cm_{\cg\setminus 1 \to 2}. 
\]
This follows by applying \cref{corollary:graphical component} to the statement $1 \indep 2 | 5$, since in each of these graphs, we can see that every d-connecting path from $1$ to $2$ contains the edge $1 \to 2$. Thus we get that
\[
\cm_\cg \cap \cm_{1 \indep 5 | 4} \cap \cm_{1 \indep 2 |5} = \cm_{\cg \setminus 1 \to 2}. 
\]

The order in which one adds the statements may change how one applies our criteria, but it obviously does not change the final result. For example if we first add the statement $1 \indep 2 | 5$, then we see that
\[
\cm_\cg  \cap \cm_{1 \indep 2 |5} = \cm_{\cg \setminus 1 \to 2}
\]
by \cref{corollary:graphical component}. We then immediately see that in  $\cg \setminus 1 \to 2$, it clearly holds that $4$ d-separates $1$ from $5$. Thus we again get that $\cm_\cg \cap \cm_{1 \indep 5 | 4} \cap \cm_{1 \indep 2 |5} = \cm_{\cg \setminus 1 \to 2}.$ 

\end{remark}

We note though that the iterative version of our technique may fail to recognize that a model $\cm_\cg \cap \cm_\mathcal{C}$ decomposes into a union of graphical models in some cases. In particular, it could be that adding the first statement $i \indep j | K$ yields a model $\cm_\cg \cap \cm_{i \indep j | K}$ which is not graphical but the addition of a later statement $a \indep b | C$ may yield a model $\cm_\cg \cap \cm_{i \indep j | K} \cap \cm_{a \indep b | C}$ which is a union of graphical models. 

\section{Implications in Gaussian Almost-Graphical CI Models}\label{sec:gaussoid}
In this section we provide a short discussion and some conjectures on the difficulty of solving \cref{prob:CIimplicationOnDAGs} which concerns our restricted class of models $\gmci$ versus the more general CI implication problem stated in \cref{prob:CIimplication}.  It is well-established that if $\mathcal{C}$ is an arbitrary set of CI statements, then to check if $a \indep_\Sigma b | C$ for all $\Sigma \in V_\mathcal{C} = \overline{\cm_\mathcal{C}}$, one must check that
$|\Sigma_{aC, bC}| = 0$ for all $\Sigma \in V_\mathcal{C}$ \cite{TobiasThesis, algstat2018}. This requires computing a \emph{Gr\"obner basis} for the ideal $\langle |\Sigma_{aC, bC}| ~:~ a \indep b | C \in \mathcal{C} \rangle$ which can be doubly-exponential in the number of variables \cite{coxlittleoshea}. In contrast, \cref{cor:CIimplicationViaPrincipal} shows that this problem can be solved on the level of algebraic varieties by evaluating determinants and then testing principal ideal membership which is significantly easier. We now introduce a conjecture which states that solving \cref{prob:CIimplicationOnDAGs} may be even easier. We begin by recalling the \emph{gaussoid} axioms. 

\begin{definition}[Gaussoid Axioms]
\label{defn:GaussoidAxioms}
Let $X \sim \mathcal{N}(\mu, \Sigma)$ be a multivariate Gaussian random vector. Then the following implications hold
\begin{enumerate}
\item $i \indep j | L \text{ and } i \indep k | jL \implies i \indep k | L \text{ and } i \indep j | kL$
\item $i \indep j | kL \text{ and } i \indep k | jL \implies i \indep j | L \text{ and } i \indep k | L$
\item $i \indep j | L \text{ and } i \indep k | L \implies i \indep j | kL \text{ and } i \indep k | jL$
\item $i \indep j | L \text{ and } i \indep j | kL \implies i \indep k | L \text{ or } j \indep k | L$
\end{enumerate}
for all distinct $i,j,k \in [n]$ and $L \subseteq [n] \setminus \{i,j,k\}$. 
\end{definition}

It is well-known that these axioms hold for all Gaussian distributions. That is, if a Gaussian distribution satisfies the premise of any of these implications, then it must satisfy the conclusion as well. Moreover, for graphical models any CI statement which holds for the model can be deduced from the set local Markov statements of $\cg$ by applying the gaussoid axioms \cite{lauritzen2018unifying}. The following conjecture suggests that the gaussoid axioms may be sufficient to determine all CI statements which are true for the model $\gmci$. 

\begin{conjecture}
\label{conj:GaussoidAxiomsSuffice}
Let $\cg$ be a DAG and $i \indep j | K$ be a conditional independence statement such that $i \not \perp_\cg j | K$. If $a \indep b | C$ holds for all $\Sigma \in \gmci$, then $a \indep b | C$ is implied by applying the gaussoid axioms to the set of CI statements $\glob(\cg) \cup \{i \indep j | K\}$. 
\end{conjecture}

\begin{example}
Let $\cg$ be the DAG pictured in \cref{figure:running-example} and consider adding the statement $2 \indep 4 | 5$ to the model $\cm_\cg$ to get the model $\cm_{\cg, 2 \indep 4 | 5}$. By applying \cref{cor:CIimplicationViaPrincipal}, we can see that $1 \indep_\Sigma 4 |5$ also holds for all $\Sigma \in \cm_{\cg, 2 \indep 4 | 5}$. However, note that $1 \perp_\cg 4 | \{2, 5\}$. Thus by the first axiom in \cref{defn:GaussoidAxioms}, we see that if $i = 4, j = 2, k = 1$ and $L = \{5\}$ then 
\[
4 \indep 2 | 5 \text{ and } 4 \indep 1 | \{2,5\}  \implies 4 \indep 1 | 5 \text{ and } 4 \indep 2 | \{1, 5\}
\]
thus we see that the new statement which we detected algebraically is already implied by applying the semigraphoid axiom. Furthermore, one can check with \cref{cor:CIimplicationViaPrincipal} that these two statements are the only additional statements which hold for $\gmci$ so in this case the semigraphoid axiom suffices to determine all additional CI statements which hold for $\gmci$. 
\end{example}

The following proposition shows that \cref{conj:GaussoidAxiomsSuffice} holds for all models involving only four random variables.

\begin{proposition}
\label{prop:ConjTrueFor4nodes}
\cref{conj:GaussoidAxiomsSuffice} holds for all models $\gmci$ where $\cg$ has four nodes. 
\end{proposition}
\begin{proof}
In \cite[Section 4]{boege2019geometry}, the authors show that up to symmetry, the only CI structures which are not closed under the gaussoid axioms for positive definite matrices are
\begin{itemize}
\item $\{1 \indep 2|3,~ 1 \indep 3|4,~ 1 \indep 4 | 2\}$,
\item $\{1\indep 2,~ 1\indep 2|\{3, 4\},~ 3 \indep 4 | 1, ~ 3 \indep 4 | 2\}$,
\item $\{1 \indep 2,~ 1 \indep 3|\{2, 4\},~ 2 \indep 4 | \{1, 3\},~  3 \indep 4\}$,
\item $\{1\indep 2|3,~ 1 \indep 3 | 4,~ 2 \indep 4 | 1,~ 3 \indep 4 | 2\}$,
\item $\{1\indep 2,~ 1 \indep 3|\{2, 4\},~ 2\indep 4 | 3,~ 3 \indep 4 | 1\}$.
\end{itemize}
This immediately implies that that if $\glob(\cg) \cup \{i \indep j | K\}$ is not equivalent to one of the previous five sets of CI statements up to symmetry, then \cref{conj:GaussoidAxiomsSuffice} holds by \cite[Theorem 4.1]{boege2019geometry}. For any of the above five sets of CI statements $\mathcal{C}$, it is impossible to find a graph $\cg$ such that the $\glob(\cg)$ is equivalent to a subset of $\mathcal{C}$ of cardinality $|\mathcal{C}| - 1$ under the gaussoid axioms. We prove this by simply checking directly for each subset of any of the potential CI structures $\mathcal{C}$, that there is no DAG $\cg$ on four nodes together with a statement $i \indep j | K$ such that $\glob(\cg) \cup \{i \indep j | K\}$ is equivalent to $\mathcal{C}$ under the gaussoid axioms. 
\end{proof}

\begin{figure}
\centering
\begin{tikzpicture}[>=stealth',shorten >=1pt,auto,node distance=0.8cm,scale=.9, transform shape,align=center,minimum size=3em]
\node[state] (w1) at (0,0) {$1$};
\node[state] (w2) at (0,3) {$2$};
\node[state] (w3) at (3,3) {$3$};
\node[state] (w4) at (3,0) {$4$};

\path[->]   (w1) edge   (w3);
\path[->]   (w1) edge   (w4);
\path[->]   (w3) edge   (w2);
\path[->]   (w2) edge   (w4);
\path[->]   (w3) edge   (w4);
\end{tikzpicture}
\caption{A graph for whose global Markov property consists only of the CI statement $1 \indep 2 | 3$.}
\label{fig:semi graphoid}
\end{figure}
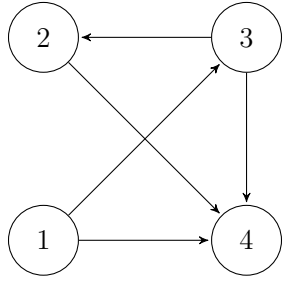

While we conjecture that the gaussoid axioms suffice to solve the CI implication problem for all models $\gmci$, it is easy to see that if one adds two additional statements to a graphical model then this immediately breaks. For instance, if $\cg$ is the graph pictured in \cref{fig:semi graphoid} then adding the statements $\{1 \indep 3 |4, 1\indep 4 |2\}$ yields the first of the CI structures in \cref{prop:ConjTrueFor4nodes}. Lastly, we note that it may be interesting to investigate \cref{prob:CIimplicationOnDAGs} and an appropriate analogue of \cref{conj:GaussoidAxiomsSuffice} in the discrete setting where fewer conditional independence inference rules will hold.

\section{Strong Faithfulness}\label{sec:approximate CI implication}

In Sections \ref{sec: unfaithful}, \ref{sec:graphicalImplication} and \ref{sec:gaussoid}, we have investigated the implication problem where we add a single conditional independence statement to a graphical model. As discussed in the introduction this can be seen as studying the algebraic geometry of a faithfulness violation. The notion of faithfulness has, however, been criticised as artificial since near conditional independence cannot be differentiated from exact conditional independence in finite samples. In response the notion of $\lambda$-faithfulness has been proposed as an alternative, which requires that for triples $i,j,K$ such that $i \not\indep_{G} j| K$, $|\rho_{i,j\cdot K}| > \lambda$ for some $\lambda >0$ as opposed to just $|\rho_{i,j\cdot K}| > 0$ \cite{zhang2002strong}. This raises the natural question: how does a graphical model to which we add one near conditional independence statement, that is, a $\lambda$-faithfulness violation behave? This problem is also interesting from an information theoretic perspective. The conditional independence constraints in a graphical model can be thought of as encoding conditional mutual information constraints of the form $I(i,j | K)=0$ but it is unclear which additional information constraints a graphical model imposes \cite{cover2012elements}. One natural class to study are constraints of the form  $I(i,j| K) -  I(i,l| K) \geq 0$. Since
    \begin{align*}
        I(i,j| K) -  I(i,l| K) 
        = \frac{1}{2} \log\left(\frac{1-\rho^2_{i,j\cdot K}}{1-\rho^2_{i,l\cdot K}}\right)
    \end{align*}
such a constraint holds precisely when any $\lambda$-faithfulness violation of the form $|\rho_{i,j\cdot K}| \leq \lambda$ implies that $|\rho_{i,l\cdot K}| \leq \lambda$, i.e., when a $\lambda$-faithfulness violation propagates \cite{arellano2013shannon}. A $\lambda$-faithfulness violation corresponds to imposing an additional inequality constraint on our model and is therefore more difficult to study with algebraic tools. Thus, we restrict ourselves to studying whether approximate implication holds for correlations of the form $\rho_{i,l\cdot K}$ only and do so with more analytical tools. The following result shows that for such correlations, a remarkably simple necessary and sufficient graphical criterion characterizing approximate implication exists.

\begin{theorem}
Let $G=(V,E)$ be a DAG and $\{i\},\{j\},K \subset V$ be sets of vertices such that $K$ does not $d$-separate $i$ from $j$.
Let $\mathcal{M}_{i \perp_{\delta} j | K}$ denote the model of covariance matrices such that $|\rho_{i,j\cdot K}|\leq \delta$, and consider a node $l$ such that  $K$ does not $d$-separate $i$ from $l$. Then for all $\delta \in [0,1]$, $\cm_\cg \cap \mathcal{M}_{i \perp_{\delta} j | K} \subseteq \mathcal{M}_{i \perp_{\delta} l | K}$ if and only if $i \indep_G l | \{j\} \cup K$.
\label{theorem: approximate}
\end{theorem}

\begin{proof}
    To allow the use of $i,j$ and $k$ as indices, let $x=i, a=j,B=K$ and $c=l$. We first show that if $c \perp_{\cg} x | a \cup B$, then $|\rho_{x,c\cdot B}| \leq |\rho_{x,a\cdot B}|$. Suppose $c \perp_{\cg} x | a \cup B$. Then, $\sigma_{x,c\cdot B}=\sigma_{x,c\cdot a,B} + \sigma_{x,a\cdot B}\sigma^{-1}_{a,a\cdot B}\sigma_{a,c\cdot B} = \sigma_{x,a\cdot B}\sigma^{-1}_{a,a\cdot B}\sigma_{a,c\cdot B}$ and therefore $$\rho_{x,c\cdot B}=\frac{\sigma_{x,c\cdot B}}{\sqrt{\sigma_{x,x\cdot B}\sigma_{c,c\cdot B}}}=\frac{\sigma_{x,a\cdot B}}{\sqrt{\sigma_{x,x\cdot B}\sigma_{a,a\cdot B}}}\frac{\sigma_{a,c\cdot B}}{\sqrt{\sigma_{a,a\cdot B}\sigma_{c,c\cdot B}}} = \rho_{x,a\cdot B} \rho_{a,c\cdot B}.$$ Since $|\rho_{a,c\cdot B}| 
\leq 1$ it follows that $|\rho_{x,c\cdot B}| \leq |\rho_{x,a\cdot B}|$.

We now show that if $c \not\perp_{\cg} x | B$ and $c \not\perp_{\cg} x | a \cup B$ hold, then we can construct a distribution $M$ such that $|\rho_{x,a\cdot B}|$ is arbitrarily small while $|\rho_{x,c\cdot B}|$ is not. 

Since $c \not\perp_{\cg} x | B$ it follows that there exists a path $p_1$ from $x$ to $c$ that we can choose by \cref{lemma: wlog} to contain colliders $n^1_1,\dots,n^1_k$ such that there exist directed paths $q^1_1,\dots,q^1_k$, each from $n^1_i$ to some $b^1_i \in B$ with possibly $n^1_i=b^1_i$, $p_1$ and the $q^1_i$s only intersect at the $n^1_i$s and the $q^1_i$s do not intersect each other. Suppose $a$ lies on neither $p_1$ nor any of the $q^1_i$s. Consider the model $\mathcal{M} \subset \cm_\cg$ where all edge coefficients on $p_1$ and the $q^1_i$s are randomly sampled, all remaining ones are set to $0$, and all error variances are randomly sampled from a distribution absolutely continuous with respect to Lebesgue. For any $\Sigma \in \mathcal{M}$, $a \indep x | B$ and for almost all $\Sigma \in \mathcal{M}$, $c \not\indep x | B$ concluding our proof. We can therefore assume that $a$ lies on either $p_1$ or some $q^1_j$ for the remainder of this proof. 
Suppose first that the latter is the case. In this case $p_1$ is also open given $a \cup B$. Consider again a covariance matrix $\Sigma \in \mathcal{M}$. Let $\Sigma(\epsilon)$ be the covariance matrix we obtain when we replace the error variances for nodes on the segment $q^1_j(a,b^1_j)$ with the value $\epsilon/k$ and the edge coefficients on the same segment are set to $1$, where $k=|q^1_j(a,b^1_j)|$ is the number of edges on the path $q^1_j(a,b^1_j)$ compared to the model that defined $\Sigma$.

In the model corresponding to $\Sigma(\epsilon)$, $b^1_j = a + z$ with $z$ a variance $\epsilon$ random variable independent of $x,a$ and $B_{-j}=B\setminus b^1_j$. Therefore, $\sigma_{x,x \cdot B_{-j}},\sigma_{a,a \cdot B_{-j}},\sigma_{c,c\cdot B_{-j}},\sigma_{x,c \cdot B_{-j}}$ do not depend on $\epsilon$. Further, $\sigma_{x,a \cdot B_{-j}}=\sigma_{x,b_j \cdot B_{-j}},\sigma_{a,b_j \cdot B_{-j}}=\sigma_{a,a \cdot B_{-j}}$ and $\sigma_{b_j,b_j \cdot B_{-j}}=\sigma_{a,a \cdot B_{-j}}+\epsilon$, where to ease notation $b^1_j=b_j$.

Therefore,
\begin{align*}
    \sigma_{x,a\cdot B} 
    &= \sigma_{x,a \cdot B_{-j}} - \sigma_{x,b_{j} \cdot b_{-j}} \sigma_{b_j,b_j \cdot B_{-j}}^{-1} \sigma_{b_j,a\cdot B_{-j}}\\
    &= \sigma_{x,a \cdot B_{-j}}  - \sigma_{x,a \cdot B_{-j}}  (\sigma_{a,a \cdot B_{-j}} + \epsilon )^{-1} \sigma_{a,a \cdot B_{-j}}.  
\end{align*}
Similarly,
\begin{align*}
    \sigma_{a,a\cdot B} 
    &= \sigma_{a,a \cdot B_{-j}} - \sigma_{a,b_j \cdot B_{-j}} \sigma_{b_j,b_j \cdot B_{-j}}^{-1} \sigma_{b_j,a\cdot B_{-j}}\\
    &= \sigma_{a,a \cdot B_{-j}} - \sigma_{a,a \cdot B_{-j}} (\sigma_{a,a \cdot B_{-j}}+\epsilon )^{-1} \sigma_{a,a \cdot B_{-j}},          
\end{align*}
while 
\begin{align*}
    \sigma_{x,x\cdot B} 
    &= \sigma_{x,x \cdot B_{-j}} - \sigma_{x,b_j \cdot B_{-j}} \sigma_{b_j,b_j \cdot B_{-j}}^{-1} \sigma_{b_j,x \cdot B_{-j}}\\
    &= \sigma_{x,x \cdot B_{-j}} - \sigma_{x,a \cdot B_{-j}} (\sigma_{a,a \cdot B_{-j}}+ \epsilon )^{-1} \sigma_{a,x \cdot B_{-j}}.       
\end{align*}
Based on this,  $\lim_{\epsilon\rightarrow 0}\sigma_{x,a\cdot B}/\sqrt{\sigma_{a,a\cdot B}}=0$, while $\lim_{\epsilon\rightarrow 0}\sigma_{x,x\cdot B}=\sigma_{x,x \cdot B_{-j}}\neq 0$. As a result, $$\lim_{\epsilon\rightarrow 0}\rho_{x,a\cdot B}=\lim_{\epsilon\rightarrow 0}\sigma_{x,a\cdot B}/\sqrt{\sigma_{x,x\cdot B}\sigma_{a,a\cdot B}}=0.$$ On the other hand, 
\begin{align*}
    \sigma_{c,c\cdot B} 
    &= \sigma_{c,c\cdot B_{-j}} - \sigma_{c,b_j\cdot B_{-j}} \sigma_{b_j,b_j \cdot B_{-j}}^{-1} \sigma_{b_j,c \cdot B_{-j}}\\
    &= \sigma_{c,c\cdot B_{-j}} - \sigma_{ca.b_{-j}} (\sigma_{a,a \cdot B_{-j}}+ \epsilon )^{-1} \sigma_{a,c \cdot B_{-j}},      
\end{align*}
and
\begin{align*}
    \sigma_{x,c\cdot B} 
    &= \sigma_{x,c \cdot B_{-j}} - \sigma_{x,b_j \cdot B_{-j}} \sigma_{b_j,b_j \cdot B_{-j}}^{-1} \sigma_{b_j,c \cdot B_{-j}}\\
    &= \sigma_{x,c \cdot B_{-j}} - \sigma_{x,a \cdot B_{-j}} (\sigma_{a,a \cdot B_{-j}}+ \epsilon )^{-1} \sigma_{a,c \cdot B_{-j}}.       
\end{align*}
and therefore $\lim_{\epsilon\rightarrow 0}\rho_{x,c\cdot B}=\lim_{\epsilon\rightarrow 0}\sigma_{x,c\cdot B}/\sqrt{\sigma_{x,x\cdot B}\sigma_{c,c\cdot B}}=\rho_{x,c.a,B_{-j}} \neq 0$. 

We can therefore assume for the remainder of this proof that $a$ lies on $p_1$. More broadly, we can assume that there does not exist a path from $x$ to $c$ that is open given $B$ and $a \cup B$. In particular, this implies that any path from $x$ to $c$ open given $B$ must contain $a$ as a non-collider, and any path open given $a \cup B$ must contain a collider $N$, such that $N \in \an(a)$ and $N \notin \an(B)$. 

Since $c \not\perp_{\cg} x | a \cup B$ there exists a path $p_2$ from $x$ to $c$ with $(l+1)$ colliders $n^2_{1},\dots,n^2_{l},n_a$ with corresponding directed paths $q^2_{1},\dots,q^2_{l},q^2_a$, ending each respectively with  $b^2_1,\dots,b^2_{l},a$ such that $p_2$ and the $q^2$s intersect at most once and the $q^2$s do not intersect. We can use these paths to construct paths $p_3=p_2(x,n_a)\oplus q_A(n_a,a)$ and $p_4=q_a(a,n_a)\oplus p_2(n_a,c)$ that are both open given $B$ and end with an edge into $a$. We will now use the paths $p_1,p_2,p_3$ and $p_4$ along along with the restriction on paths that are open given $B$ to show that $\cg$ must contain a specific structure. We will then exploit this structure to construct a submodel $\mathcal{M} \in \cm_\cg$ with the required properties.

Let $f^1_i$ denote the $i$-th fork node on $p_1$, i.e., a node of the form $\leftarrow f^1_i \rightarrow$ of which there are $k+1$. Since every node on the subpath $p_1(f^1_1,f^1_{k+1})$ is an ancestor of $B$, $a$ must either lie on the segment i) $p_1(x,f^1_1)$ or ii) $p_1(f^1_{k+1},c)$.

Case i): Let $i_2$ be the node closest to $c$ on $p_1$, where $p_1$ and $p_3$ intersect and consider the path $p_5=p_3(x,i_2) \oplus p_1(i_2,c)$. If $i_2 \in B$, then $i_2\neq a$ and $i_2$ must be a collider on both $p_1$ and $p_3$ and therefore on $p_5$. As a result, $p_5$ is open given $B$ and $a \cup B$. We can therefore, without loss of generality, assume that $i_2 \notin B$. If $i_2 \neq a$ and $i_2 \notin B$, then $p_5$ does not contain $a$ and therefore must be closed given $B$, which requires that $i_2 \notin \an(B)$. Since all nodes on $p_1(f^1_1,f^1_{k+1})$ are ancestors of $B$ by construction, it follows that $i_2$ lie on $p_1(x,f_1)$. If $i_2 = a$, then it is by nature of $p_1$ and $p_3$ a collider on $p_5$ and we again obtain that $a$ must lie on $p_1(x,f_1)$. Furthermore, we can replace $p_2$ with $p_5$ since $p_5$ is a path from $x$ to $c$ that is open given $a \cup B$ but closed given $B$, i.e., we can assume $p_1$ and $p_2$ intersect at some node $i_2$ such that they have the segment $p_1(i_2,c)$ in common. 

Let $i_2' \neq i_2$ be the node closest to $i_2$ on $p_1(a,i_2)$, where $p_1$ and $p_2$ intersect. Then we can again replace $p_2$ with $p_2(x,i'_2) \oplus p_1(i'_2,c)$. We can therefore without loss of generality assume that there is no such $i_2'$, i.e., $p_1$ and $p_2$ do not intersect between $a$ and $i_2$ except at $i_2$. Let $i_1\neq i_2$ be the node closest to $i_2$ on $p_2(x,i_2)$ where $p_1(x,i_2)$ and $p_2(x,i_2)$ intersect and consider $p=p_1(x,i_1) \oplus p_2(i_1,c)$. Since $i_1$ lies on the directed towards $x$ path $p_1(x,i_2)$ it must be a non-collider on $p$ and since $p_1(x,i_2)$ does not contain any nodes in $B$ it follows that $p$ is open given $a \cup B$. We can therefore replace $p_2$ with $p$. In summary we can assume that there exists nodes $i_1$ and $i_2$ such that $p_1(x,i_1)=p_2(x,i_1)$, $p_1(i_2,c)=p_2(i_2,c)$ and $p_1$ and $p_2(i_1,i_2)$ only intersect at $i_1$ and $i_2$. It may be the case that $i_1=x$ or $i_2= a$. \cref{fig:i} illustrates the resulting structure.

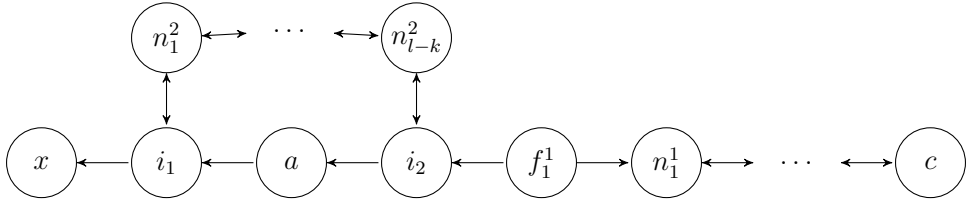
\begin{figure}[t!]
	\centering
		\begin{tikzpicture}[>=stealth',shorten >=1pt,auto,node distance=0.8cm,scale=.9, transform shape,align=center,minimum size=3em]
		\node[state] (w1) at (0,0) {$x$};
        \node[state] (w2) [right =of w1] {$i_1$};
		\node[state] (w3) [right =of w2] {$a$}; 
		\node[state] (w4) [right =of w3] {$i_2$};
		\node[state] (w5) [right =of w4] {$f^1_1$};
		\node[state] (w6) [right =of w5] {$n^1_1$};
		\node[] (w7) [right =of w6] {$\dots$};
		\node[state] (w8) [right =of w7] {$c$};
  
  \node[state] (t1) [above =of w2] {$n^2_{1}$};
  \node[] (t2) [above =of w3] {$\dots$};
  \node[state] (t3) [above =of w4, label=center:$n^2_{l-k}$] {};


  \path[<-]   (w1) edge    (w2);
\path[<-]   (w2) edge    (w3);
\path[<-]   (w3) edge    (w4);
\path[<-]   (w4) edge    (w5);
\path[->]   (w5) edge    (w6);
\path[<->]   (w6) edge    (w7);
\path[<->]   (w7) edge    (w8);

\path[<->]   (t1) edge    (w2);
\path[<->]   (t3) edge    (w4);
\path[<->]   (t1) edge    (t2);
\path[<->]   (t2) edge    (t3);


		\end{tikzpicture}
	\caption{Figure illustrating the Case i). The bi-directed edges represent paths of the form $l \leftarrow \dots \leftarrow f \rightarrow \dots \rightarrow  r$.}
	\label{fig:i}
\end{figure}

To further characterize the possible structures we will now show that we can without loss of generality also restrict how any $q^1_i$ path intersects with $p_2$, any $q^2_j$ path intersects with $p_1$ and finally, how any $q^1_i$ path intersects any $q^2_j$ path.
Consider a path $q^1_i$ and suppose it intersects with $p_2(i_1,i_2)$ and let $i'$ be the node closest to $i_1$ on $p_2(i_1,i_2)$ where this is the case. Then, $p=p_2(x,i') \oplus q^1_i(i',n^1_i) \oplus p_1(n^1_i,c)$ is open given $B$ and does not contain $a$, where we use that $B \cap \de(i') \neq \emptyset$, that $n^1_i$ is a non-collider on this path and finally, that if $i' 
\in B$ it is a collider on $p$. We can therefore, without loss of generality, assume that the paths $q^1_i,i \in \{1,\dots,k\}$ do not intersect with $p_2(i_1,i_2)$.

Consider a path $q^2_i,i \in \{1,\dots,l-k\}$ and suppose it intersects with $p_1(i_1,i_2)$. Let $i'$ be the node closest to $i_1$ on $p_1(i_1,i_2)$ where this is the case. Then $p=p_1(x,i') \oplus q^2_{i}(i',n^2_i) \oplus p_2(n^2_i,c)$ is open given $B$ and does not contain $a$, where we use that $B \cap \de(i') \neq \emptyset$, that $n^2_i$ is a non-collider on this path and finally, that if $i' 
\in B$ it is a collider on $p$. We can therefore, without loss of generality, assume that the $q^2_i$ for $i \in \{1,\dots,l-k\}$ do not intersect with $p_1$.

Finally suppose for some $i \in \{1,\dots,k\}$ and some $j \in \{1,\dots,k-l\}$, $q^1_i$ and $q^2_j$ intersect and let $i$ be the node closest to $n^2_j$ on $q^2_j$ where this is the case. Then $p_2(x,n^2_j) \oplus q_j(n^2_j,i) \oplus q_i(i,n^1_i) \oplus p_1(n^1_i,c)$ is a path from $x$ to $c$ that does not contain $a$ and is open given $B$, since $i \in \an(B)$ is a collider on $p$.  We can therefore, without loss of generality, assume that $q^1_i$ does not intersect with $q^2_j$. In summary, $p_1,p_2,$ the $q^1_i$s and the $q^2_j$s (for $j \in \{1,\dots,l-k\}$) only intersect at the $n_i$'s.

We will now construct a linear structural equation model.
Consider the linear structural equation model $M$ where all edge coefficients on $p_1,p_2$, the $q^1_i$s and the $q^2_j$s are $1$, except those adjacent to the $f^1_i$s and the $f^2_j$s which are alternatingly $1$ and $-1$ beginning with $1$ and all remaining ones are set to $0$. Further assume that all error variances for $i_1,i_2,C,f^1_1,\dots,f^1_{k},f^2_{1},\dots,f^2_{l-k}$ are $1$ and all remaining ones are $0$, i.e., $x=i_1,a=i_2$ and $v^k_i=n^k_i=-f^k_i + f^k_{i+1}$ with $k=1,2$. Let $M(\epsilon)$ be the same model except the error variances for $f^2_{1},\dots,f^2_{l-k}$ are set to $\epsilon$. Let $\tilde{G}$ be the corresponding simplified graph with which $M(\epsilon)$ is compatible.

Let $B'=\{b^1_1,\dots,b^1_k\}$ and $B'' = \{b^2_1,\dots,b^2_{l-k}\}$.
Since $C \perp_{\tilde{\cg}} B'' | B',\sigma_{B'',c\cdot B'}=0$. As a result, 
\begin{align*}
    \sigma_{x,c\cdot B}
    &=\sigma_{x,c\cdot B'}-\sigma_{x,B''\cdot B'}\Sigma_{B'',B''\cdot B'}^{-1} \sigma_{B'',c\cdot B'} \\ 
    &=\sigma_{x,c\cdot B'}
\end{align*}
does not depend on $\epsilon$ and similarly for $\sigma_{c,c\cdot B}$. 

The terms $\sigma_{x,x\cdot B},\sigma_{a,a\cdot B}$ and $\sigma_{x,a\cdot B}$ do depend on $\epsilon$.
By construction, however $\sum_{i=1}^{k-l+1} b^2_{i} = - f^2_1 +f^2_{k-l+1}$ and $\sum_{i=1}^{k} b^1_{i} = - f^1_1 + f^1_{k+1}$. Therefore, if $\beta=(-1,\dots,-1)$, then 
\begin{align*}
     x-\beta B
     &= \epsilon_{i_1} + \epsilon_{i_2} + f^2_{1} -f^2_{k-l+1} + f^1_{1} - f^2_{1} +f^2_{k-l+1} - f^1_1 + f^1_{k+1} \\
     &= \epsilon_{i_1} + \epsilon_{i_2} + f^1_{k+1}.
\end{align*}
Therefore, the variance of the residual $x-\beta B$ does not depend on $\epsilon$ if and only if $\beta=(-1,\dots,-1)$ and as a result, the ordinary least squares regression coefficient $\beta_{x,B}$ (which minimizes the variance of $X - \beta_{x,B}B$) has to converge to $\beta$ as $\epsilon$ goes to infinity. 
Therefore, $\lim_{\epsilon \rightarrow \infty}\sigma_{x,x\cdot B} = 3$. Finally, using that  $\sigma_{x,a} = \sigma_{a,a}= Var(\epsilon_{i_1}) + Var(f^2_{k-l+1}) + Var(f_1^1) = 2 + \epsilon$ and $\sigma_{B,a}=(0,\dots,-\epsilon,-1,\dots,0)$, it follows that $\lim_{\epsilon \rightarrow \infty} \sigma_{x,a\cdot B} = \lim_{\epsilon \rightarrow \infty} (\sigma_{x,a} - \beta_{x,B} \sigma_{B,a})=1$. Regarding $\sigma_{a,a\cdot B}$, there exists no $\beta$ such that $a-\beta B$ does not depend on $f^2_{1}, \dots, f^2_{l-k}$ and therefore $\sigma_{a,a\cdot B}$ goes to infinity as $\epsilon$ goes to infinity. Combined, this implies that by choosing a large $\epsilon$ we can make $\rho_{x,a\cdot B}=\sigma_{x,a\cdot B}/\sqrt{\sigma_{x,x\cdot B}\sigma_{a,a\cdot B}}$ arbitrarily small, while 
$\rho_{x,c\cdot B}=\sigma_{x,c\cdot B}/\sqrt{\sigma_{x,x\cdot B}\sigma_{c,c\cdot B}}\geq c > 0$. This concludes Case i).

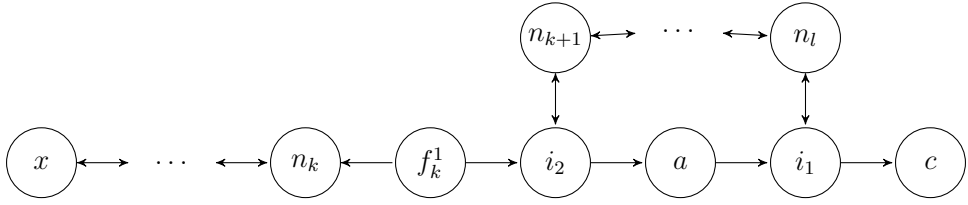
\begin{figure}[t!]
	\centering
		\begin{tikzpicture}[>=stealth',shorten >=1pt,auto,node distance=0.8cm,scale=.9, transform shape,align=center,minimum size=3em]
		\node[state] (w1) at (0,0) {$x$};
        \node[] (w2) [right =of w1] {$\dots$};
		\node[state] (w3) [right =of w2] {$n_k$}; 
		\node[state] (w4) [right =of w3] {$f^1_k$};
		\node[state] (w5) [right =of w4] {$i_2$};
		\node[state] (w6) [right =of w5] {$a$};
		\node[state] (w7) [right =of w6] {$i_1$};
		\node[state] (w8) [right =of w7] {$c$};
  
  \node[state] (t1) [above =of w5, label=center:$n_{k+1}$] {};
  \node[] (t2) [above =of w6] {$\dots$};
  \node[state] (t3) [above =of w7] {$n_{l}$};


  \path[<->]   (w1) edge    (w2);
\path[<->]   (w2) edge    (w3);
\path[<-]   (w3) edge    (w4);
\path[->]   (w4) edge    (w5);
\path[->]   (w5) edge    (w6);
\path[->]   (w6) edge    (w7);
\path[->]   (w7) edge    (w8);

\path[<->]   (t1) edge    (w5);
\path[<->]   (t3) edge    (w7);
\path[<->]   (t1) edge    (t2);
\path[<->]   (t2) edge    (t3);


		\end{tikzpicture}
	\caption{Figure illustrating the Case ii). The bi-directed edges represents paths of the form $l \leftarrow \dots \leftarrow f \rightarrow \dots \rightarrow  r$.}
	\label{fig:ii}
\end{figure}

Case ii): Using $p_4$ instead of $p_3$ we can argue as in Case i) to arrive at the structure illustrated in \cref{fig:ii}. From there we can again proceed as in Case i) to construct a model class $M(\epsilon)$ such that for $\epsilon \rightarrow \infty, \sigma_{a,a\cdot B} \rightarrow \infty$ while $\sigma_{x,x\cdot B},\sigma_{c,c\cdot B},\sigma_{x,c\cdot B},\sigma_{x,a\cdot B}$ converge to non-zero constants.
\end{proof}

\begin{lemma}
Let $G=(V,E)$ be a DAG and $\{x\},\{y\},K \subset V$ be sets of vertices such that such that $x\perp_{\cg}y| K$. Then there exists a path $p$ from $x$ to $y$ such that i) no non-collider is in $K$, ii) for every collider $n_i$ on $p$ there exists a causal path $q_i$ from $n_i$ to some $k_i \in K$, iii) each $q_i$ intersects with $p$ only at $n_i$ and iv) the $q_i$s do not intersect with each other.
    \label{lemma: wlog}
\end{lemma}

\begin{proof}
    By the definition of d-separation, i) and ii) are trivial. Regarding iii) see the proof of Lemma B.3 in the Appendix of \cite{henckel2022graphical}. Regarding iv), suppose that $q_i$ and $q_j$ intersect where without loss of generality we assume that $n_i$ is to the right of $n_j$ on $p$. Let $i$ be the node closest to $n_i$ at which $q_i$ and $q_j$ intersect. Then $p(x,n_i) \oplus q(n_i,i) \oplus q_j(i,n_j) \oplus p(n_j \oplus Y)$ is path from $x$ to $y$ that is open given $K$ with at least one less collider than $p$ ($i$ is a collider but $n_i$ and $n_j$ are non-colliders on $p$). Repeatedly applying this procedure we either obtain a colliderless path from $x$ to $y$, or a path with colliders whose $q'$s respect condition iv).
\end{proof}

\begin{lemma}
    Consider a path $p_1$ from $a$ to $b$ and path $p_2$ from $c$ to $d$ in a DAG $\cg$ that are both open given some set $K$ and intersect. Let $i$ be the node closest to $a$ on $p_1$ where the two intersect and consider the path $p=p_1(a,i)\oplus p_2(i,d)$. Then $p$ is closed given $K$ if and only if $K \cap \de(i) = \emptyset$ and $i$ is a collider on $p$. 
    \label{lemma: aux open}
\end{lemma}

\begin{proof}
    By choice if $i$, $p_1(a,i)$ and $p_3$ only intersect at $i$ and therefore $p$ is a path. If $i=a$ or $i=d$, then $p$ is a subpath of $p_1$, respectively $p_2$ and therefore trivially open given $K$. If $i \in K$, then it must be a collider on both $p_1$ and $p_2$ and therefore also on $p$ and as a result $p$ is open given $K$. If $i \notin Z$ and $i$ is a non-collider on $p$, then $p$ is trivially open given $Z$ and similarly, if $i$  is a collider and $K \cap \de(i) \neq \emptyset$. As a result, $p$ may only be closed given $K$ if $i$ is a collider on $p$ and $K \cap \de(i) = \emptyset$.
\end{proof}

\begin{example}
    Let $\cg$ be the DAG pictured in \cref{figure:running-example} and consider adding the statement $2 \indep 4$ to the model $\cm_\cg$ to get the model $\cm_{\cg, 2 \indep 4}$. By \cref{proposition: monomial} and \cref{theorem:d-connection monomial}, $\cm_{\cg, 2 \indep 4}$ decomposes into a union of graphical models, specifically the model corresponding to the graph $\cg$ without the edge $1 \rightarrow 2$ and the graph without the edge $2 \rightarrow 4$. It is easy to verify that in both graphs $1$ is d-separated from $5$ and therefore, $1 \indep_\Sigma 5$ for all $\Sigma \in \cm_{\cg, 2 \indep 4}$. However, $1$ is not d-separated from $5$ given $4$ and therefore given a $\Sigma \in \cm$, $|\rho_{14}| \leq \lambda$ does not imply that $|\rho_{15}| \leq \lambda$ by \cref{theorem: approximate}. This is an example where a faithfulness violation propagates and a $\lambda-$faithlessness violation does not, which illustrates that the two are distinct.
\end{example}

\section{Discussion}
In this paper we introduce a new restricted version of the Gaussian conditional independence implication problem which only concerns sets of conditional independence statements of the form $\glob(\cg) \cup \{i \indep j | K\}$. We study the exact version of this problem as well as the approximate version where conditional independence is replaced with small correlations. 

In the exact case, we show that this problem is essentially equivalent to a principal ideal membership problem but requires computing determinants of matrices with large polynomial entries which can be exponential in the size of the matrix. However, we give a complete characterization of when the model $\gmci$ decomposes as a union of graphical models and is thus easily solvable with d-separation. In the approximate case, we show that for a certain class of statements implication holds if and only if a simple d-separation statement holds. This indicates that this problem may be significantly easier than the traditional implication problem.

Lastly, we note that our work may shed light on how certain causal discovery algorithms perform when presented with additional CI statements which violate the typical assumption of faithfulness to some DAG. In particular, the methods developed in \cref{sec:graphicalImplication} provide a first step toward the problem of determining if a collection of CI statements can be faithfully encoded by a graphical model (or a union of graphical models).

\section*{Acknowledgements}
The authors are grateful for helpful discussions with Tobias Boege and Seth Sullivant. Mathias Drton and Pratik Misra received funding from the European Research Council (ERC) under the European Union’s Horizon 2020 research and innovation programme (grant agreement No. 883818). Benjamin Hollering was supported by the Alexander von Humboldt Foundation.

\bibliographystyle{plain}
\bibliography{mybib.bib}

\newpage

\appendix

\section{Additional preliminaries}
\label{appendix: prelims}

\textbf{Graphical notation:} A simple directed graph $\cg=(V,E)$ consists of nodes $V$ and directed edges $E$, i.e., edges of the form $i\rightarrow j$ such that there is at most one edge between any two nodes. We say that two nodes $i$ and $j$ are adjacent if they are connected by an edge. A path $p$ is a sequence if distinct nodes $(i, \dots, k)$ such that for all consecutive nodes on $p$ are adjacent. We say $p$ is directed if all edges on $p$ are of the form $i \rightarrow j$, i.e, point towards $k$. A directed path from $i$ to $k$ and the edge $k \rightarrow i$ form a directed cycle. We call a simple directed graph without directed cycles a directed acyclic graph (DAG). Given a path $p=(i,\dots j,\dots,l,\dots,k)$, let $p(j,l)$ denote the subpath from $j$ to $l$. Given two paths $p=(i,\dots,k)$ and $q=(k,\dots,l)$ let $p \oplus q=(i,\dots,k,\dots,l)$ denote the concatanation of the two paths. Given an edge $i\rightarrow j$ we say that $i$ is a parent of $j$ and denote the set of parents with $\pa(j)$. Given a directed path $i \rightarrow \dots \rightarrow j$ we say that $i$ is an ancestor of $j$, respectively $j$ is a descendant of $i$. We denote the set of ancestors with $\an(j)$ and the set of descendants with $\de(i)$.

\textbf{Covariance matrix notation:} Consider mutually Gaussian random vectors $A,B$ and $C$. We denote the covariance matrix between $A$ and $B$ with $\Sigma_{A,B}$ and the conditional covariance of $A$ and $B$ given $C$, i.e. $\Sigma_{A,B}-\Sigma_{A,C} \Sigma_{C,C}^{-1} \Sigma_{C,B}$,  with $\Sigma_{A,B\cdot C}$. For singleton $a$ or $b$ we write $\sigma_{a,b\cdot C}$. Similarly, we denote the conditional correlation between $A$ and $B$ given $C$, i.e. $\Sigma_{A,A\cdot C}^{-1/2} \Sigma_{A,B \cdot C} \Sigma_{B,B\cdot C}^{-1/2}$, with $\rho_{A,B \cdot C}$. For the ordinarily least squares coefficient of the regression of $A$ on $B$ we write $\beta_{A,B}=\Sigma_{B,B}^{-1} \Sigma_{B,A}$.

\end{document}